\documentstyle[amsfonts,amscd,epsfig]{amsart}

\numberwithin{equation}{section}

\newtheorem{theorem}{Theorem}[section]

\newtheorem{lemma}[theorem]{Lemma}
\newtheorem{proposition}[theorem]{Proposition}

\newtheorem{definition}{Definition}[section]

\newcommand{\sn}{S_n}
\newcommand{\snme}{S_{n-1}}
\newcommand{\snmeme}{S_{n-2}}

\newcommand{\mnme}{{\cal M}_{0,n-1}}

\newcommand{\mn}{{\cal M}_{0,n}}
\newcommand{\hnme}{H^*(\mnme, \Bbb{Q})}

\newcommand{\hn}{H^*(\mn, \Bbb{Q})}
\newcommand{\hni}{H^i(\mn, \Bbb{Q})}

\begin{document}
\title{A formula for the Euler characteristic of $\overline{{\cal M}}_{2,n}$}
\author{G. Bini}
\author{G. Gaiffi}
\author{M. Polito}

\address{Scuola Normale Superiore\\
P.zza dei Cavalieri, 7\\
56126 Pisa, Italia
 \bigskip}

\email[Gilberto Bini]{bini@@cibs.sns.it}
\email[Giovanni Gaiffi]{gaiffi@@cibs.sns.it}
\email[Marzia Polito]{polito@@cibs.sns.it}

\begin{abstract}
In this paper we compute the generating function for the Euler characteristic of the 
Deligne-Mumford compactification of the moduli space of smooth $n$-pointed
genus $2$ curves. The proof relies on quite elementary methods, such as the enumeration of the graphs
involved in a suitable stratification of $\overline{{\cal M}}_{2,n}$.
\end{abstract}

\maketitle
\section{Introduction} 
\label{intro}

Let${\cal \ }\overline{{\cal M}}_{g,n}$ , $2g+n-2>0$, denote the
Deligne-Mumford compactification of the moduli space of smooth $n$-pointed
genus $g$ curves. As explained in \cite{AC}, there exists a stratification 
\begin{equation*}
\overline{{\cal M}}_{g,n}\supseteq \partial \overline{{\cal M}}%
_{g,n}\supseteq \partial ^2\overline{{\cal M}}_{g,n}\supseteq \ldots
\supseteq \partial ^{3g-3+n}\overline{{\cal M}}_{g,n}\text{,}
\end{equation*}
whose codimension $k$ strata $\partial ^k\overline{{\cal M}}_{g,n}$, $1\leq
k\leq 3g-3+n$, can be thoroughly described by genus $g$ graphs with a finite set $V$
of vertices, a finite set of edges $E$ and characterised in the following
way:

\begin{itemize}
\item  each vertex can have labelled half-edges, which we shall refer to as
leaves, with labels in the set $\left\{ 1,\ldots ,n\right\} $;

\item  there exists a map $\gamma :V\rightarrow \left\{ 0,\ldots ,g\right\} $;

\item  for each vertex $v$, $2 \gamma (v)+h(v)+l(v)-2>0$, where $l(v)$ is the
number of edges issuing from $v$ and $h(v)$ is the number of leaves stemming
from $v$;

\item $g= \sum_{v \in V} \gamma (v) +$ the topological genus of the graph.
\end{itemize}

\noindent Observe that a graph with $\gamma (v)=0$, for every $v$, and without loops
is simply a tree.

Let us now consider $X_j:=\partial ^j\overline{{\cal M}}_{g,n}-\partial
^{j+1}\overline{{\cal M}}_{g,n}$, $0\leq j\leq 3g-4+n$, ($\partial ^0%
\overline{{\cal M}}_{g,n}:=\overline{{\cal M}}_{g,n}$), which is easily seen
to be a union of quasi-projective subvarieties 
$\{ X_{i} \} _{i \in I \left( j \right)}$.
With each $X_i$ we can
associate a graph which describes all the elements contained in a
codimension $j$ stratum, but not in its closure. We will make use of this
stratification to calculate generating functions for the Euler
characteristic of $\overline{{\cal M}}_{1,n}$, $\overline{{\cal M}}_{2,n}$.
More precisely, given a graph $\Gamma _i$ associated with $X_i$ there exists
a morphism 
\begin{equation*}
\xi _{\Gamma _i}:\Pi _{v\in V_{\Gamma _i}}\overline{{\cal M}}%
_{g(v),h(v)+l(v)}\rightarrow \overline{{\cal M}}_{g,n}
\end{equation*}

\noindent which is described as follows. Each vertex $v$ of the graph
corresponds to a smooth $(h(v)+l(v))$- pointed genus $g(v)$ curve $\left[
C;x_1,\ldots ,x_{h(v)},x_{h(v)+1},\ldots ,x_{l(v)}\right] $ which is
attached to another curve, corresponding to another vertex, through one of
the marked points $x_j$, $h(v)+1\leq j\leq l(v)$. By the multiplicativity
of the Euler characteristic, this means that $\chi (X_i)$ is generically given by the
product of Euler characteristics of each moduli space corresponding to
vertices of $\Gamma _i$. In some cases, depending on the symmetries of the graph,
we should take into account a symmetric group action.

By the additivity of the Euler characteristic, we just have to find out how
to enumerate all the graphs involved in the complete description of boundary
strata. To this end, we introduce the generating function for trees, which
satisfies the following recursive relation: 
\begin{equation*}
D(t):=t+\sum_{n\geq 2}\chi ({\cal M}_{0,n+1})\frac{D\left( t\right) ^n}{n!}%
\text{,}
\end{equation*}

\noindent and observe that our graphs are obtained by some particular
configurations to which trees are attached (see sections \ref{graphs} 
and \ref{k2} for these
configurations and the definition of {\bf graph-type}).

We can finally state our main results. Consider the generating function 
\begin{equation*}
K_g(t):=\sum_{2g+n-2>0}\chi \left( \overline{{\cal M}}_{g,n}\right) \frac{t^n%
}{n!}\text{.}
\end{equation*}

\noindent Then the following theorems hold.

\begin{theorem}
\label{genus1} 
\begin{eqnarray*}
K_1\left( t\right) \ =\frac{19}{12}D+\frac{23}{24}D^2+\frac 5{18}D^3+\frac{%
D^4}{24}-\frac E{12}-\frac 12\log \left( 1-\log \left( 1+D\right) \right) 
\text{.}
\end{eqnarray*}
\end{theorem}

\begin{theorem}
\label{genus2} 
\begin{eqnarray*}
\ K_2 \left( t \right) &=&\frac 1{1440\left( 1+D\right) ^2\left( E-1\right) ^3}
[-2D^8 {\left( E - 1 \right) } {\left( 7+3E \right) 
-24D^7\left( E-1\right) ^2\left( 17E - 7 \right) }
\\
& &-3D^6 \left( E-1 \right)^2 \left( 201E+259 \right) 
+30D^5 \left( E-1 \right) ^2 \left( 61E-221 \right) \\
& &+15D^4\left( 631E^3-2640E^2+3395E-1386\right) \\
&&+60D^3\left(341E^3-1322E^2+1633E-652\right)\\
&&+180D^2\left( 138E^3-519E^2+635E-254\right) \\
&&+360D\left( 45E^3-167E^2+206E-84\right)   \\
&& +60\left(73E^3-270E^2+336E-144\right) ] \text{,}\\
\end{eqnarray*}
where we have set $E:=\log \left( 1+D\right) $.
\end{theorem}

Although a generating function for $\chi (\overline{{\cal M}}_{1,n})$ has
already been found by Getzler in \cite{G1}, we propose an elementary and
simplified algorithm to determine it. This will allow the reader to become
more familiar with calculations for the genus $2$ case. While preparing this
work, we came to know that a generating function for $\chi \left( \overline{%
{\cal M}}_{g,n}\right) $ has been found by J. Harer; even if our computation
works in theory for every genus $g$, at the moment it turns out to be
effective only for small genera, but in these cases it provides very direct
algebraic formulas for the sought-for generating functions. We also notice
that for the computation we need to know the characteristic of the open sets 
${\cal M}_{g,n}\subseteq \overline{{\cal M}}_{g,n}$ (when $g=0,1,2$):
although a general method to compute $\chi ({\cal M}_{g,n})$ can be found in 
\cite{HZ}, our paper is self-contained, since we shall make explicit
calculations for $\chi ({\cal M}_{0,n})$, $\chi ({\cal M}_{1,n})$ and $\chi (%
{\cal M}_{2,n})$ on the basis of a convenient stratification as suggested in 
\cite{AC}. We also checked that our results, for low $n$,
coincide with the ones obtained by E. Getzler in \cite{G2}, and with some recent
computations of him, which he kindly informed us about.

We would like to thank E. Arbarello, M. Cornalba and J. Harer for their help and 
encouragement during the development of our work.

Many of our explicit computations were performed using 
{\it Mathematica}\footnote{\copyright {\it Mathematica} is a registered trademark of
Wolfram Research, Inc.}; 
we are grateful to D. Finocchiaro for his patience in teaching us to use it.

\section{Quotients of products of ${\cal M}_{0,n} $ and ${\cal M}_{1,n} $}

\label {quoz}

As remarked in the Introduction, we shall compute $\chi ($ $\overline{{\cal M%
}}_{1,n})$ and $\chi (\overline{{\cal M}}_{2,n})$ summing the contributions
provided by the Euler characteristic of the subvarieties $X_i$ represented
by graphs. When computing the characteristic associated with each graph, we
have to deal with the characteristic of ${\cal M}_{0,n}$ and some of its
quotients with respect to the symmetric group action which permutes the
marked points. Furthermore, we will need quotients of products of ${\cal M}%
_{0,n}$ with respect to suitable actions of the symmetric group. In
addition, these quotients are also useful to describe strata in the
stratification introduced in sections \ref{m2n} and \ref{m1n} in order to compute 
directly the characteristic of the open sets ${\cal M}_{1,n}$ and ${\cal M}_{2,n}$.

This section is devoted to these computations, which are of independent
interest and will be carried out using two different techniques. On one
side, one can study the $S_n$-invariants in the rational cohomology ring of $%
{\cal M}_{0,n}$ (subsection \ref{invariant}), and on the other side one can employ
geometric arguments involving branched covering (subsection \ref{geometric}).
In this last subsection, we are going to use some results of section \ref{m1n}, namely
the formula for the Euler characteristic of ${\cal M}_{1,n}$; nevertheless this formula
will be proved without assuming any of the previous results of this paper.

Obviously, we start by recalling the Euler characteristic of ${\cal M}_{0,n}$.
This is provided by the following

\begin{theorem}
For $n\geq 3$%
\begin{equation*}
\chi ({\cal M}_{0,n})=(-1)^{n-3}(n-3)!\text{.}
\end{equation*}
\end{theorem}

{\bf Proof.} Consider the fibration $\pi :{\cal M}_{0,n+1}\rightarrow {\cal M%
}_{0,n}$ with fiber ${\Bbb P}^1-\{n$ points$\}$. This gives the recursive
formula 
\begin{equation*}
\chi ({\cal M}_{0,n+1})=(2-n)\chi ({\cal M}_{0,n})\text{,}
\end{equation*}

\noindent with initial data $\chi ({\cal M}_{0,3})=1$.

\begin{flushright}
$\Box $
\end{flushright}

Let us now consider the case of the symmetric group $S_j$ acting on ${\cal M}%
_{0,n}$ (here, for $j\leq n$, we identify $S_j$ with a subgroup of $S_n$).
We notice that when $n-j\geq 3$ the quotient map $q:{\cal \ M}%
_{0,n}\rightarrow {\cal M}_{0,n}/S_j$ is unramified, since any automorphism
of ${\Bbb P}^1$ fixing three or more points is the identity. This implies
that 
\begin{equation*}
\chi ({\cal M}_{0,n}/S_j)=\frac{\chi ({\cal M}_{0,n})}{j!}=\frac{%
(-1)^{n-3}(n-3)!}{j!}\text{.}
\end{equation*}

The following subsection provides a description of the $S_n$ module $H^{*}(%
{\cal M}_{0,n},{\Bbb Q})$ which allows us to compute the Euler
characteristic of the quotient spaces ${\cal M}_{0,n}/S_n$, ${\cal M}%
_{0,n}/S_{n-1}$, ${\cal M}_{0,n}/S_{n-2}$.

\medskip\ 

\subsection{$S_n$-invariants of $H^{*}({\cal M}_{0,n},{\Bbb Q})$}
\label{invariant}

The following well known theorem of Invariant Theory points out the
relations between the cohomology ring of $\mn/S_j$ ($1\leq j \leq
n$) and the $S_j$-invariants in $\hn$.
\begin{theorem}
Let $X$ be a variety and $G$ a finite group which acts
 on $X$.
Then

\begin{displaymath}
H^* ( X/G, \, {\Bbb Q})  \cong (H^*(X,\, {\Bbb Q}))^G.
\end{displaymath}

\label{inva}
\end{theorem}
\noindent Let us first recall some results about the symmetric group
action
on $\hn$.
For every $n\geq 3$ and every $j\leq n$ we will denote  by $Ch_j(\mn)$
(resp. $Ch_j^i(\mn)$)
the character of the $S_j$
representation on $\hn$ (resp. on $\hni$).
Furthermore, we will denote by $I_j$ and $P_j$ respectively the characters
of the trivial and  standard representations of
$S_j$. Then we have
\begin{theorem}{(see \cite{Ga}, \cite{Ma})}
For every $n\geq 3$,
\begin{equation}
        Ch_{n-1}^i(\mn)=Ch_{n-1}^i(\mnme)+P_{n-1} Ch_{n-1}^{i-1}(\mnme).
\label{matieu}
\end{equation}
\end{theorem}
We will also refer to the following theorem which was first obtained by
Lehrer in \cite{L} and that can be derived by formula  (\ref{matieu}).
\begin{theorem}
For every $n\geq 3$,
\begin{displaymath}
        Ch_{n-1}(\mn)=Ind^{\snme}_{S_2}(I_2).
\end{displaymath}
\label{indotto}
\end{theorem}

Let us now denote by $(\, , \,)_{S_j}$ the inner  product
in the space of class functions on $S_j$ (in the sequel we may omit the
subscript $S_j$ if it is clear to which group we are referring to).

\begin{lemma}
For $n\geq 3$,
\begin{equation}
        (Ch_n(\mn), I_n)_{\sn}=1
        \label{e1}
\end{equation}
\begin{equation}
        (Ch_{n-1}(\mn), I_{n-1})_{\snme}=1
        \label{e2}
\end{equation}
\begin{equation}
        (Ch_{n-2}(\mn), I_{n-2})_{\snmeme}=n-2.
        \label{e3}
\end{equation}

\label{almeno}
\end{lemma}
{\bf Proof}. This is a consequence of Theorem \ref{indotto}.
In fact we can write
\begin{displaymath}
        (Ch_{n-1}(\mn), I_{n-1})_{\snme}=(Ind^{\snme}_{S_2}(I_2),
        I_{n-1})_{\snme}
\end{displaymath}
which, by Frobenius reciprocity law, is equal to
\begin{displaymath}
        (I_2, Res^{\snme}_{S_2}(I_{n-1}))_{S_2}=(I_2,
I_2)_{S_2}=1
\end{displaymath}
This gives relation (\ref{e2}).
As for relation (\ref{e1}) we note that,  since
$Res^{\sn}_{\snme}(Ch_{n}(\mn))$ is equal to
$Ch_{n-1}(\mn)$, then
\begin{displaymath}
        1=(Ch_{n-1}(\mn), I_{n-1})_{\snme}=
\end{displaymath}
\begin{displaymath}
         = (Ch_{n}(\mn),
        Ind^{\sn}_{\snme}(I_{n-1}))_{\sn}= (Ch_{n}(\mn),
        I_n + P_n)_{\sn}
\end{displaymath}
Since dim $H^0(\mn, {\Bbb Q})=1$ we know that $(Ch_{n}(\mn),
        I_n )_{\sn}\geq 1$. This implies that $(Ch_{n}(\mn),
        I_n )_{\sn}=1$ and $(Ch_{n}(\mn),
         P_n)_{\sn}=0$.

It remains to prove the last assertion, which can be formulated as
$$(Res^{\snme}_{\snmeme}(Ch_{n-1}(\mn)),
        I_{n-2})_{\snmeme}=n-2.$$
Then we can write
\begin{displaymath}
        (Res^{\snme}_{\snmeme}(Ch_{n-1}(\mn)),
        I_{n-2})_{\snmeme}=(Ch_{n-1}(\mn),
        I_{n-1}+P_{n-1})_{\snme}=
\end{displaymath}
\begin{displaymath}
        = 1+ (Ch_{n-1}(\mn),
        P_{n-1})_{\snme}
\end{displaymath}
which, applying Theorem \ref{indotto}, is equal to
\begin{displaymath}
        =1+ (I_2, Res^{\snme}_{S_2}(P_{n-1}))_{S_2}.
\end{displaymath}
The second addendum can be easily computed using the branching rule and
is equal to $n-3$. This completes the proof.
\begin{flushright}
$\Box $
\end{flushright}

We are now ready to compute Euler characteristics.
\begin{theorem}
The Euler characteristic of $\mn/S_n$ and of $\mn/\snme$ is equal
to 1 for every $n\geq 3$.

The Euler characteristic of $\mn/S_{n-2}$ is equal to 0 if $n$ is
even and equal to 1 is $n$ is odd.

\end{theorem}
{\bf Proof}. The assertions concerning $\mn/S_n$ and  $\mn/\snme$
are immediate
consequences of Theorem \ref{inva} and Lemma \ref{almeno}, since 
the one dimensional spaces
$H^0(\mn/S_n,\Bbb{Q})$ and $H^0(\mn/\snme,\Bbb{Q})$ afford respectively
 the trivial representations $I_n$ and $I_{n-1}$.

Turning to  $\mn/S_{n-2}$, we want to prove by induction on $n$ the
following stronger proposition which implies our claim:

\begin{proposition}
 For every $n\geq 3$ and $0\leq i \leq n-3$, in the irreducible
decomposition  of the $\snmeme$ module $H^i(\mn,\Bbb{Q})$ the trivial
representation $I_{n-2}$ occurs exactly with multiplicity 1.
\end{proposition}
 {\bf Proof.} The base of  induction ($n=3$) is obvious. Now, given $n>3$, it
suffices to prove that, for every $i$, the multiplicity of $I_{n-2}$ in
the decomposition of $H^i(\mn,\Bbb{Q})$ is at least 1. In fact we  then
observe that the top cohomology of $\mn$ has degree $n-3$ and we can
apply Lemma \ref{almeno}. For $i=0$  our assertion is trivial. Let us then suppose $i\geq 1$.
From the relation of Theorem \ref{matieu} we deduce
 \begin{displaymath}
 Ch_{n-2}^i(\mn)=
        Ch_{n-2}^i(\mnme)+      Res^{\snme}_{\snmeme}(P_{n-1})
                Ch_{n-2}^{i-1}(\mnme),
 \end{displaymath}
 that is to say,
 \begin{equation}
  Ch_{n-2}^i(\mn)=
        Ch_{n-2}^i(\mnme)+      (P_{n-2}+I_{n-2})
                Ch_{n-2}^{i-1}(\mnme)
                \label{utile}
 \end{equation}
 If $i= 1$, we have $Ch_{n-2}^{i-1}(\mnme)=I_{n-2}$, therefore $I_{n-2}$
 appears in the irreducible decomposition of  $Ch_{n-2}^i(\mn)$.
 If $i\geq 2$, we observe that, by  the inductive hypothesis,
 $(Ch_{n-3}^{i-1}(\mnme), I_{n-3})=1$.
 But by Frobenius reciprocity law
 \begin{displaymath}
        1=(Ch_{n-3}^{i-1}(\mnme), I_{n-3})=(Ch_{n-2}^{i-1}(\mnme), I_{n-2}
+
        P_{n-2}).
 \end{displaymath}
 Now $(Ch_{n-2}^{i-1}(\mnme), I_{n-2})=0$ since we have already proven
 that the only subspace of $\hnme$ which affords the trivial
 representation $I_{n-2}$ is $H^0(\mnme,\Bbb{Q})$. Then we have
 \begin{displaymath}
        1=(Ch_{n-2}^{i-1}(\mnme),
        P_{n-2}).
 \end{displaymath}
 Therefore, in the equation (\ref{utile}) we find  the
 product $P_{n-2}P_{n-2}$ as an addendum: its  decomposition into
irreducibles (see
 \cite{FH} Chap. 4) is equal to $I_{n-2} + P_{n-2}$ plus two other
 irreducible characters. $\Box $
 
\smallskip\

\subsection{The geometric method}
\label{geometric}

Let us now apply the fibration method to other quotients of ${\cal M}_{0,n}$.
 A first case is given by the action of the Klein group $S_2\times S_2$
which intertwines two pairs of marked points, i.e.

\begin{equation*}
\left[ {\Bbb P}^1;x_1,\ldots ,x_{n-1},x_n\right] \rightarrow \left[ 
{\Bbb P}^1;x_1,\ldots ,x_n, x_{n-1}\right] \text{.}
\end{equation*}

\begin{proposition}
\label{esseduedue}For $n\geq 4$%
\begin{equation*}
\chi ({\cal M}_{0,n}/(S_2\times S_2))=\left\{ 
\begin{array}{l}
0\text{, }n=4\text{,} \\ 
\\ 
0\text{, }n=5\text{,} \\ 
\\ 
-2\text{, }n=6\text{,} \\ 
\\ 
\frac{\chi ({\cal{M}}_{0,n})}4\text{, }n\geq 7\text{.}
\end{array}
\right.
\end{equation*}
\end{proposition}

{\bf Proof. }For $n\geq 7$, the claim follows easily, since there are at
least three fixed points. The remaining cases can be treated in the same way
and we show $\chi ({\cal M}_{0,6}/(S_2\times S_2))$ as a sample case. The
quotient map $\varphi :{\cal M}_{0,6}\rightarrow {\cal M}_{0,6}/(S_2\times
S_2)$ is a $4-1$ covering ramified over the branch locus 
\begin{equation*}
{\cal B}:=\left\{ \left[ {\Bbb P}^1;0,\infty ,1,-1,b,-b\right] \text{; }%
b\neq 0,\infty ,1,-1\right\} \text{.}
\end{equation*}

In fact, the fiber of a generic point $\left[ {\Bbb P}^1;0,\infty
,1,a,b,c\right] $ of ${\cal M}_{0,6}/(S_2\times S_2)$ contains the points

$\left[ {\Bbb P}^1;0,\infty ,1,a,b,c\right] $, $\left[ {\Bbb P}^1;0,\infty
,1,a,c,b\right] $, $\left[ {\Bbb P}^1;0,\infty ,1,1/a,b/a,c/a\right] $, and

\noindent $\left[ {\Bbb P}^1;0,\infty ,1,1/a,c/a,b/a\right] $. These are distinct
points except when $a=-1$, $b=-c$ in which case the fiber is made up of two
points. We now have the following relation 
\begin{equation*}
\chi ({\cal M}_{0,6})=4\chi ({\cal M}_{0,6}/(S_2\times S_2))-2\chi ({\cal B})%
\text{.}
\end{equation*}
Since the Euler characteristic of the branch locus is $-1$ ( as being
isomorphic to $\left( {\Bbb P}^1-\{0,\infty ,1,-1\}\right) /S_2$), we deduce
that $\chi ({\cal M}_{0,6}/(S_2\times S_2))=-2$.

\begin{flushright}
$\Box $
\end{flushright}

Another important example is given by the action of the dihedral group $D_4$
on ${\cal M}_{0,n}$. Here we identify $D_4$ with the subgroup of $S_4$
generated by $\sigma $ and $\tau $ such that 
\begin{equation*}
\sigma \cdot \left[ {\Bbb P}^1;x_1,x_2,x_3,x_4,\ldots ,x_n\right]
\rightarrow \left[ {\Bbb P}^1;x_2,x_1,x_3,x_4,\ldots ,x_n\right] \text{,}
\end{equation*}

\begin{equation*}
\tau \cdot \left[ {\Bbb P}^1;x_1,x_2,x_3,x_4,\ldots ,x_n\right] \rightarrow
\left[ {\Bbb P}^1;x_3,x_4,x_1,x_2,\ldots ,x_n\right] \text{.}
\end{equation*}

\begin{proposition}
For $n\geq 4$%
\begin{equation*}
\chi ({\cal M}_{0,n}/D_4)=\left\{ 
\begin{array}{l}
0\text{, }n=4\text{,} \\ 
\\ 
0\text{, }n=5\text{,} \\ 
\\ 
-1\text{, }n=6\text{,} \\ 
\\ 
\frac{\chi ({\cal{M}}_{0,n})}8\text{, }n\geq 7\text{.}
\end{array}
\right.
\end{equation*}
\end{proposition}

{\bf Proof. }The proof proceeds analogously as in Proposition \ref
{esseduedue}. In the case $n=6$, the branch locus is exactly the same and so
the result follows.

\begin{flushright}
$\Box $
\end{flushright}

Let us now come to a detailed analysis of those quotients of products we
shall use. Take $n_1$, $n_2\geq 3$ positive integers such that $n_1\leq n_2$ and
consider the product ${\cal M}_{0,n_1}\times {\cal M}_{0,n_2}$ endowed with
the action of $S_2$ which permutes a pair of marked points in each element
of both factors of the product.

\begin{proposition}
\label{prod1} 
\begin{equation*}
\chi (({\cal M}_{0,n_1}\times {\cal M}_{0,n_2})/S_2)=\left\{ 
\begin{array}{l}
\chi (({\cal M}_{0,n_2})/S_2)\text{, }n_1=3\text{,} \\ 
\\ 
\chi (({\cal M}_{0,n_2})/S_2)-\chi ({\cal M}_{0,n_2})\text{, }n_1=4\text{,}
\\ 
\\ 
\frac{\chi ({\cal M}_{0,n_1})\chi ({\cal M}_{0,n_2})}2\text{, }n_1\geq 5%
\text{.}
\end{array}
\right.
\end{equation*}
\end{proposition}

{\bf Proof.} Since ${\cal M}_{0,3}$ is a point, the first statement follows
easily. On the other hand, when $n_1\geq 5$, the Euler characteristic of the
product is simply given by the relation 
\begin{equation*}
\chi ({\cal M}_{0,n_1}\times {\cal M}_{0,n_2})=2\chi (({\cal M}%
_{0,n_1}\times {\cal M}_{0,n_2})/S_2)\text{,}
\end{equation*}

\noindent since in this case the projection of $({\cal M}_{0,n_1}\times 
{\cal M}_{0,n_2})$ onto the quotient does not have any branch points. For
the second case, we recall that ${\cal M}_{0,4}$ maps onto ${\cal M}%
_{0,4}/S_2$ with a branch point, namely $[{\Bbb P}^1;0,\infty ,1,-1]$.
Therefore, because of the action of $S_2$ on ${\cal M}_{0,4}\times {\cal M}%
_{0,n_2}$, we have 
\begin{equation*}
\chi ({\cal M}_{0,4}\times {\cal M}_{0,n_2})=2\chi (({\cal M}_{0,4}\times 
{\cal M}_{0,n_2})/S_2)-\chi (U_{n_2})\text{,}
\end{equation*}

\noindent where $U_{n_2}$ is the branch locus of the map ${\cal M}%
_{0,n_2}\rightarrow {\cal M}_{0,n_2}/S_2$. Since 
\begin{equation*}
\chi ({\cal M}_{0,n_2})=2\chi (({\cal M}_{0,n_2})/S_2)-\chi (U_{n_2})\text{,}
\end{equation*}

\noindent we finally deduce that 
\begin{equation*}
\chi (({\cal M}_{0,4}\times {\cal M}_{0,n_2})/S_2)=\chi (({\cal M}%
_{0,n_2})/S_2)-\chi ({\cal M}_{0,n_2})\text{.}
\end{equation*}

\begin{flushright}
$\Box $
\end{flushright}

Let us consider the action of $S_3$ on the product ${\cal M}_{0,n_1}\times 
{\cal M}_{0,n_2}$ by permuting triples of marked points in each element of
the product. Then we can prove the following

\begin{lemma}
i) $\chi (({\cal M}_{0,3}\times {\cal M}_{0,n_2})/S_3)=\chi (({\cal M}%
_{0,n_2}/S_3))$,

ii) $\chi (({\cal M}_{0,4}\times {\cal M}_{0,4})/S_3)=2$, $\chi (({\cal M}%
_{0,5}\times {\cal M}_{0,4})/S_3)=1$, $\chi (({\cal M}_{0,5}\times {\cal M}%
_{0,5})/S_3)=2$,

$\,$iii) $\chi (({\cal M}_{0,n_1}\times {\cal M}_{0,n_2})/S_3)=\frac
16\,\chi ({\cal M}_{0,n_1})\chi ({\cal M}_{0,n_2})$, when at least one of
the $n_j^{\text{'}}$s is greater or equal to $6$.
\end{lemma}

{\bf Proof.} Statements i) and iii) are obvious.

For the proof of ii), let us examine directly what happens when $%
(n_1,n_2)=(4,4)$, $(n_1,n_2)=(4,5)$, and $(n_1,n_2)=(5,5)$. We first observe
that the map 
\begin{equation*}
{\cal M}_{0,4}\rightarrow {\cal M}_{0,4}/S_3\text{,}
\end{equation*}

\noindent is ramified along the fiber with two points $p_1=\left[ {\Bbb P}%
^1;0,\infty ,1,\alpha ,\alpha ^2\right] $ and $p_2=\left[ {\Bbb P}%
^1;0,\infty ,1,\alpha ^2,\alpha \right] $, ($\alpha $ a primitive third root
of unity), and the fiber with three points $q_1=\left[ {\Bbb P}^1;0,\infty
,1,-1\right] $, $q_2=\left[ {\Bbb P}^1;0,\infty ,1,1/2\right] $, $q_3=\left[ 
{\Bbb P}^1;0,\infty ,1,2\right] $. This allows to compute directly branch
points of the map 
\begin{equation*}
{\cal M}_{0,4}\times {\cal M}_{0,4}\rightarrow ({\cal M}_{0,4}\times {\cal M}%
_{0,4})/S_3\text{,}
\end{equation*}

\noindent yielding that 
\begin{equation*}
\chi (({\cal M}_{0,4}\times {\cal M}_{0,4}))=6\chi (({\cal M}_{0,4}\times 
{\cal M}_{0,4})/S_2)-11\text{,}
\end{equation*}

\noindent i.e. $\chi (({\cal M}_{0,4}\times {\cal M}_{0,4})/S_2)=2$.

\noindent For the other cases, we recall that the map 
\begin{equation*}
{\cal M}_{0,5}\rightarrow {\cal M}_{0,5}/S_3
\end{equation*}

\noindent has a fiber with two points, instead of six, namely $v_1=\left[ 
{\Bbb P}^1;0,\infty ,1,\alpha ,\alpha ^2\right] $ and $v_2=\left[ {\Bbb P}%
^1;0,\infty ,1,\alpha ^2,\alpha \right] $. So when $(n_1,n_2)=(4,5)$ branch
points are given by $(p_1,v_1)$ , $(p_2,v_2)$, $(p_1,v_2)$, $(p_2,v_1)$;
hence it is easy to compute $\chi (({\cal M}_{0,4}\times {\cal M}%
_{0,5})/S_3)=2$. In the same way, one can treat the remaining case.
Notice that these results could also have been obtained via representation theory of the symmetric group.

\begin{flushright}
$\Box $
\end{flushright}

\medskip\ 

We next consider the action on ${\cal M}_{0,n_1}\times {\cal M}_{0,n_2}$
 of the group $S_2\times S_2$  generated by the following involutions: 
\begin{equation*}
\sigma_{1} \left( \left[ {\Bbb P}^1;x_1,\ldots ,x_{n_1-1},x_{n_1}\right] ,\left[ {\Bbb P%
}^1;x_1,\ldots ,x_{n_2-3},x_{n_2-2},x_{n_2-1},x_{n_2}\right] \right) \rightarrow
\end{equation*}
\begin{equation*}
\rightarrow \left( \left[ {\Bbb P}^1;x_1,\ldots ,x_{n_1},x_{n_1-1}\right]
,\left[ {\Bbb P}^1;x_1,\ldots ,x_{n_2-2},x_{n_2-3},x_{n_2-1},x_{n_2}\right] \right) \text{,}
\end{equation*}

\begin{equation*}
\sigma_{2} \left( \left[ {\Bbb P}^1;x_1,\ldots ,x_{n_1-1},x_{n_1}\right] ,\left[ {\Bbb P%
}^1;x_1,\ldots ,x_{n_2-3},x_{n_2-2},x_{n_2-1}x_{n_2}\right] \right) \rightarrow
\end{equation*}
\begin{equation*}
\rightarrow \left( \left[ {\Bbb P}^1;x_1,\ldots ,x_{n_1-1},x_{n_1}\right]
,\left[ {\Bbb P}^1;x_1,\ldots ,x_{n_2-3},x_{n_2-2},x_{n_2},x_{n_2-1}\right] \right) \text{.}
\end{equation*}

\begin{proposition}
Take integers $n_1,n_2$ such that $3\leq n_1 < n_2$. Then

i) $\chi (({\cal M}_{0,n_1}\times {\cal M}_{0,n_2})/(S_2\times S_2))=\chi (%
{\cal M}_{0,n_2}/(S_2\times S_2))$, for $n_1=3$,

ii) 
\begin{equation*}
\chi (({\cal M}_{0,4}\times {\cal M}_{0,n_2})/(S_2\times S_2))=\left\{ 
\begin{array}{l}
0\text{, }n_2=4\text{,} \\ 
\\ 
-1\text{, }n_2=5\text{,} \\ 
\\ 
1\text{, }n_2=6\text{,} \\ 
\end{array}
\right.
\end{equation*}

iii) $\chi (({\cal M}_{0,n_1}\times {\cal M}_{0,n_2})/(S_2\times S_2))=\frac
12\chi ({\cal M}_{0,n_1}) \chi (({\cal M}_{0,n_2})/S_2)$, for $n_1 \geq 5$, and $n_2 \geq 4$.
\end{proposition}

{\bf Proof. }i) is obvious. For ii) one needs to recall that the only
ramification point in the quotient map ${\cal M}_{0,4}\rightarrow {\cal M}%
_{0,4}/(S_2\times S_2)$ is $\left[ {\Bbb P}^1;0,\infty ,1,-1\right] $ . In
the other cases, ramification points are given by maps introduced in
Proposition \ref{esseduedue}. Finally, iii) follows by the fact that a
rational point with more than three marked points is automorphism free.

\begin{flushright}
$\Box $
\end{flushright}

\medskip\ 

Eventually, let us now consider the group $S_2\times S_2$  
acting on ${\cal M}%
_{0,n_1}\times {\cal M}_{0,n_2}\times {\cal M}_{0,n_3}$, generated by:

\begin{equation*}
(\left[ {\Bbb P}^1;x_1,\ldots ,x_{n_1-1},x_{n_1}\right] ,\left[ {\Bbb P}%
^1;x_1,\ldots ,x_{n_2-1},x_{n_2}\right] ,\left[ {\Bbb P}%
^1;x_1,\ldots ,x_{n_3-3},x_{n_3-2},x_{n_3-1},x_{n_3}\right] )\rightarrow
\end{equation*}
\begin{equation*}
\rightarrow (\left[ {\Bbb P}^1;x_1,\ldots ,x_{n_1},x_{n_1-1}\right] ,\left[ 
{\Bbb P}^1;x_1,\ldots ,x_{n_2-1},x_{n_2}\right] ,\left[ 
{\Bbb P}^1;x_1,\ldots x_{n_3-2},x_{n_3-3},x_{n_3-1},x_{n_3}\right] )\text{,}
\end{equation*}

\begin{equation*}
(\left[ {\Bbb P}^1;x_1,\ldots ,x_{n_1-1},x_{n_1}\right] ,\left[ {\Bbb P}%
^1;x_1,\ldots ,x_{n_2-1},x_{n_2}\right] ,\left[ {\Bbb P}%
^1;x_1,\ldots ,x_{n_3-3},x_{n_3-2},x_{n_3-1},x_{n_3}\right] )\rightarrow
\end{equation*}
\begin{equation*}
\rightarrow (\left[ {\Bbb P}^1;x_1,\ldots ,x_{n_1-1},x_{n_1}\right] ,\left[ 
{\Bbb P}^1;x_1,\ldots ,x_{n_2},x_{n_2-1}\right] ,\left[ 
{\Bbb P}^1;x_1,\ldots x_{n_3-3},x_{n_3-2},x_{n_3},x_{n_3-1}\right] )\text{,}
\end{equation*}

\noindent with $n_1\geq 3$, $n_3\geq 4$, $n_2\geq 3$, $n_1\leq n_2$. In this
case we have

\begin{proposition}
i) $\chi (({\cal M}_{0,n_1}\times {\cal M}_{0,n_2}\times {\cal M}%
_{0,n_3})/(S_2\times S_2))=\chi (({\cal M}_{0,n_2}\times {\cal M}%
_{0,n_3})/(S_2\times S_2))$, for $n_1=3$.

ii) 
\begin{equation*}
\chi (({\cal M}_{0,4}\times {\cal M}_{0,4}\times {\cal M}_{0,n_3})/(S_2%
\times S_2))=\left\{ \ 
\begin{array}{l}
-1\text{, }n_3=4\text{,} \\ 
\\ 
0\text{, }n_3=5\text{,} \\ 
\\ 
-2\text{, }n_3=6\text{,} \\ 
\\ 
\frac 14\chi ({\cal M}_{0,n_3})\text{, }n_3\geq 7\text{.}
\end{array}
\right.
\end{equation*}

iii) $\chi (({\cal M}_{0,4}\times {\cal M}_{0,n_2}\times {\cal M}%
_{0,n_3})/(S_2\times S_2))=\frac 12\chi (({\cal M}_{0,4}\times {\cal M}%
_{0,n_3})/S_2)\chi ({\cal M}_{0,n_2})$, for

$n_2\geq 5$,

iv) $\chi (({\cal M}_{0,n_1}\times {\cal M}_{0,n_2}\times {\cal M}%
_{0,n_3})/(S_2\times S_2))=\frac 12\chi (({\cal M}_{0,n_2}\times {\cal M}%
_{0,n_3})/S_2)\chi ({\cal M}_{0,n_1})$, for

$n_1\geq 5$.
\end{proposition}

{\bf Proof.} The only statement to prove is ii). Let us work out the case $%
n_3=5$ as a sample case. The quotient map ${\cal M}_{0,4}\times {\cal M}%
_{0,4}\times {\cal M}_{0,5}\overset{\varphi }{\rightarrow }({\cal M}%
_{0,4}\times {\cal M}_{0,4}\times {\cal M}_{0,5})/(S_2\times S_2)$ is a
degree four map with branch locus the image under $\varphi $ of the
following set of points 
\begin{equation*}
\left\{ \left( \left[ {\Bbb P}^1;\infty ,0,1,-1\right] ,\left[ {\Bbb P}%
^1;\infty ,0,1,b,1-b\right] ,\left[ {\Bbb P}^1;\infty ,0,1,-1\right] \right)
;b\neq 0,\infty ,1,1/2\right\} \text{.}
\end{equation*}

\noindent Therefore 
\begin{equation*}
2=\chi ({\cal M}_{0,4}\times {\cal M}_{0,4}\times {\cal M}_{0,5})=4\chi ((%
{\cal M}_{0,4}\times {\cal M}_{0,4}\times {\cal M}_{0,5})/S_2\times S_2)+2%
\end{equation*}

\noindent and the claim follows.

\begin{flushright}
$\Box $
\end{flushright}

Let us consider ${\cal M}_{1,n}$, $n\geq 2$ and let the group $S_2$ act on
it by permuting the last two marked points in each genus $1$ $n$-pointed
curve.

\begin{proposition}
i) 
\begin{equation*}
\chi ({\cal M}_{1,n}/S_2)=\left\{ 
\begin{array}{l}
1\text{, }n=2\text{,} \\ 
\\ 
1\text{, }n=3\text{,} \\ 
\\ 
1\text{, }n=4\text{,} \\ 
\\ 
0\text{, }n=5\text{,} \\ 
\\ 
6\text{, }n=6\text{.}
\end{array}
\right. 
\end{equation*}

ii) $\chi ({\cal M}_{1,n}/S_2)=\frac 12\chi ({\cal M}_{1,n})$, when $n\geq 7$.
\end{proposition}

{\bf Proof.} Since an $n$-pointed genus $1$ curve has non-trivial
automorphisms when $n\geq 5$, ii) follows easily. On the other hand, each
case in i) must be analyzed separately. For $n=2$, the action of $S_2$ is
free, since the two pointed elliptic curves $\left[ C;0,p\right] $ and $%
\left[ C;p,0\right] $ , $p\neq 0$, are the same because of the structure of $%
C$ as ${\Bbb C}$ modulo a lattice.

We now turn to the case $n=3$. Suppose $\left[ C;0\right] $ is elliptic and $%
p,q$ nonzero distinct points on $C$. Clearly, if $\left[ C;0\right] \,$ is
general, the $-1$ involution around zero exchanges $p$ and $q$. Now write $C=%
{\Bbb C}/\Lambda $, and let $p$ and $q$ be the classes of the two complex
numbers $z$ and $w$. Suppose first $\Lambda ={\Bbb Z}+i{\Bbb Z}$. If the
automorphism given by multiplication by $i$ interchanges $p$ and $q$, then $%
iz\equiv w$ mod $\Lambda $ and $iw\equiv z$ mod $\Lambda $. Thus $-z\equiv z$
mod $\Lambda $, i.e. $p$ is a $2$-torsion point. Since $p$ and $q$ are
distinct, they must be the classes of $1/2$ and $i/2$. Suppose next that $%
\Lambda ={\Bbb Z}+\omega {\Bbb Z}$, where $\omega $ is a primitive third
root of unity, and let $\varphi $ be an order $6$ automorphism. If $\varphi $
interchanges $p$ and $q$, then $-p=\varphi ^3(p)=q$. In conclusion, the
branch locus of 
\begin{equation*}
{\cal M}_{1,3}\rightarrow {\cal M}_{1,3}/S_2\text{,}
\end{equation*}

\noindent consists of an isolated point plus the set $U$ of isomorphism
classes of curves $\left[ C;0,p,-p\right] $ such that $p$ is not a torsion $%
2 $-point. Since $U$ is isomorphic to ${\cal M}_{0,5}/S_3$, we finally deduce that 
\begin{equation*}
0=\chi ({\cal M}_{1,3})=2\chi ({\cal M}_{1,3}/S_2)-2\text{.}
\end{equation*}

For the other cases, we proceed analogously. When $n=4$, the branch locus of
the map 
\begin{equation*}
{\cal M}_{1,4}\rightarrow {\cal M}_{1,4}/S_2\text{,}
\end{equation*}

\noindent consists of the isolated point $\left[ C;0,i/2,1/2,1/2+i/2\right] $
and the set of points $U^{\prime }$ of isomorphism classes of curves $\left[
C;0,v,p,-p,\right] $, where $v$ is a torsion $2$-point and $p$ is not a
torsion $2$-point. Since this stratum is isomorphic to ${\cal{M}}_{0,5}/S_2$, we have 
\begin{equation*}
0=\chi ({\cal M}_{1,4})=2\chi ({\cal M}_{1,4}/S_2)-1-1\text{.}
\end{equation*}

Other computations are similar and can be worked out by the reader.

\begin{flushright}
$\Box $
\end{flushright}

\medskip\ 

We end this subsection with a mixed product. Let $S_2$ act on ${\cal M}%
_{1,n_1}\times {\cal M}_{0,n_2}$, $n_1\geq 1$, $n_2\geq 3$, as follows: 
\begin{equation*}
\left( \left[ C;x_1,\ldots ,x_{n_1-1},x_{n_1}\right] ,\left[ {\Bbb P}%
^1;y_1,\ldots ,y_{n_2-1},y_{n_2}\right] \right) \rightarrow 
\end{equation*}
\begin{equation*}
\rightarrow \left( \left[
C;x_1,\ldots ,x_{n_1},x_{n_1-1}\right] ,\left[ {\Bbb P}^1;y_1,\ldots
,y_{n_2},y_{n_2-1}\right] \right) \text{,}
\end{equation*}

\noindent where $C$ is an elliptic curve.

Similarly to Proposition \ref{prod1} one can prove: 

\begin{proposition}
\begin{equation*}
\chi (({\cal M}_{1,n_1}\times {\cal M}_{0,n_2})/S_2)=\left\{ 
\begin{array}{l}
\chi ({\cal M}_{1,n_1}/S_2)\text{, }n_2=3\text{,} \\ 
\\ 
\chi ({\cal M}_{1,n_1}/S_2)-\chi ({\cal M}_{1,n_1})\text{, }n_2=4\text{,} \\ 
\\ 
\frac 12\chi ({\cal M}_{0,n_2})\chi ({\cal M}_{1,n_2})\text{, }n_2\geq 5%
\text{.}
\end{array}
\right.
\end{equation*}
\end{proposition}

\section{The Euler characteristic of ${\cal M}_{2,n}$}
\label{m2n}

For every $n$-pointed curve $C$ of genus $2$ , we denote by $\tau $ the
hyperelliptic involution. We can stratify ${\cal M}_{2,n}$ according to the
action of $\tau $ on the marked points. In fact, we are going to decompose $%
{\cal M}_{2,n}$ into a disjoint union of quasi projective subvarieties,
which will be doubly indexed, up to isomorphism: the first index counts the
number of $\tau $-fixed marked points and the other one counts the couples
of points in involution.

In order to fix notation, we say that ${\cal M}_{2,n}$ is the union of $%
a_{j,r}$ subvarieties of type $\left\{ j,r\right\} $, isomorphic to $U_{j,r}$%
, hence 
\begin{equation*}
\chi \left( {\cal M}_{2,n}\right) =\sum_{j=0}^{\min (n,6)}\sum_{r=0}^{\left[ 
\frac{n-j}2\right] }a_{j,r}\chi \left( U_{j,r}\right) \text{.}
\end{equation*}
Set 
\begin{equation*}
U_{j,r}:=\left\{ 
\begin{array}{c}
\left[ C,p_1,...,p_n\right] :\tau (p_i)=p_i\text{ , }i=1,...,j\text{ , } \\ 
\tau (p_{j+2i})=p_{j+2i-1},i=1,...,r\text{ , }\tau (p_l)\neq p_k\text{
otherwise}
\end{array}
\right\} \text{.}
\end{equation*}

Since the subvarieties isomorphic to $U_{j,r}$ are obtained by permuting the
markings, it is easy to see by combinatorial arguments, that
\begin{equation*}
a_{j,r}=\binom nj\frac{(n-j)!}{2^r(n-j-2r)!r!}\text{.}
\end{equation*}

\noindent Let us consider the covering map 
\begin{equation*}
f_{j,r}:U_{j,r}\rightarrow \frac{{\cal M}_{0,n+6-r-j}}{S_{6-j}}
\end{equation*}
sending the class $\left[ C,p_1,...,p_n\right] $ to the class 

\begin{equation*}
\left[ C/\tau ,\left[p_1\right],...,\left[ p_j \right] ,\left[
p_{j+1}\right] = \left[ p_{j+2} \right] ,...,\left[
p_{j+2r-1} \right] = \left[ p_{j+2r} \right] ,q_1,...,q_{6-j}\right] \text{%
,}
\end{equation*}
where $\left\{ q_1,...,q_{6-j}\right\} $ are the ramification points of $%
\tau $ other than the images of the marked points, and the group $S_{6-j}$
acts permuting exactly these points.

Conversely, given such data, a $n$ pointed genus $2$ curve is determined up
to the choice of the branch of the covering where the marked points are, and
we can say that the map $f_{j,r\text{ }}$ is a covering of
degree $2^{n-j-r-1}$, except for the cases where $j=n$, and $r=0$, when it
is an isomorphism.

We claim that this map is unramified unless $j=0$, and 
$n-r=2$, where the target space is $\frac{{\cal M}_{0,8}}{S_{6}}$.
A ramification point of this map implies the existence of an isomorphism 
(different from the hyperelliptic involution $\sigma $) between two genus $2$, 
$n$-pointed curves, which represent the same equivalence 
class modulo $\sigma $; this induces an automorphism of the rational curve 
$C / \tau$, fixing the classes 
of the markings, and permuting the classes of the $\tau$-fixed points.
Thus one can see that if $j>0$ or $n-r>2$, this automorphism is the identity,
and the isomorphism between the  genus $2$ curves is the identity too, since it
fixes at least $6$ points and is not the hyperelliptic involution.

But if $j=0$ and $n-r=2$, our claim is that the ramification locus of $f_{j,r}$
is isomorphic to $\frac{{\cal M}_{0,5}}{S_{3}}$; observe that there
are only three such cases, namely $U_{0,0} \subset {\cal M}_{2,2}$, 
$U_{0,1} \subset {\cal M}_{2,3}$, $U_{0,2} \subset {\cal M}_{2,4}$. 
Let us explain the simpler case, 
namely $f_{0,0}:U_{0,0}\rightarrow \frac{{\cal M}_{0,8}}{S_{6}}$;
we need to find out for which curves $(C, p_{1}, p_{2}) \in U_{0,0}$ 
there exists an isomorphism 
$$
\sigma : (C, p_{1}, p_{2}) \rightarrow (C, p_{1}, \tau (p_{2})) \text{.}
$$

By Riemann-Hurwitz formula, and uniqueness of the hyperelliptic involution,
the curve $C / \sigma $ is elliptic, and the quotient map is ramified  exactly
over the images of $p_{1}$ and $\tau (p_{1})$. Observe that the two automorphism
$\sigma$ and $\tau$ commute, and $\tau$ (resp. $\sigma$) induces an automorphism
of $C / \sigma $ (resp. $C / \tau $). Everything  matches in the following commutative
diagram:

\begin{equation*}
\begin{CD}
C @>\phi >> C/\tau \\
@V\psi VV @V \overline{\psi} VV \\
C/\sigma @>\overline{\phi} >> C/ < \sigma, \tau > 
\end{CD}
\end{equation*}

By the conditions on $\sigma$ and $\tau$, the map $\overline{\phi}$ is ramified over 
$\overline{\psi}\phi \left( p_{2} \right)$, and over the images of the ramification points
 of $\phi$; since the map $\overline{\psi}$ has degree $2$, these six points 
 form exactly three fibers of it; the last point we have to mark on $C / < \sigma , \tau > $ 
 is $\overline{\psi}\phi \left( p_{1} \right)$; hence a genus $2$, $2$-pointed curve
satisfying our requirements determines a genus $0$ curve, with $5$ marked points, 
three of which are indistinguishable. Conversely, given a point in
 $\frac{{\cal M}_{0,5}}{S_{3}}$,  a
 point $\left[ C, p_{1}, p_{2} \right] \in U_{0,0} \subset {\cal M}_{2,2}$
 matching our conditions is uniquely determined by 
 building  $C / \tau$ and $C/ \sigma$ as ramified coverings over the marked points 
 of the rational curve, namely, by reversing the construction. This proves our claim.

Using the results of section \ref{quoz}, we have that 
\begin{equation*}
\chi \left( \frac{{\cal M}_{0,n+6-r-j}}{S_{6-j}}\right) =\left\{ 
\begin{array}{c}
\left( -1\right) ^{n+3-j-r}\frac{\left( n+3-j-r\right) !}{\left( 6-j\right) !%
}\text{ , for }n-r\geq 3\text{,} \\ 
0\text{ , for }n-r=2\text{ , }j\text{ even,} \\ 
1\text{ , for }n-r=2\text{ , }j\text{ odd,} \\ 
1\text{, for }n-r=1\text{,} \\ 
1\text{, for }n-r=0\text{.}
\end{array}
\right.
\end{equation*}

\noindent If $n\geq 7$, then the fiber of the universal curve 
\begin{equation*}
{\cal M}_{2,n+1}\rightarrow {\cal M}_{2,n}
\end{equation*}
is a genus $2$ curve without $n$ points, hence has Euler characteristic $%
-(n+2)$;for $n\geq 8$, 
\begin{eqnarray*}
\chi \left( {\cal M}_{2,n}\right) &=&\prod_{h=7}^{n-1}-(2+h)\chi \left( 
{\cal M}_{2,7}\right) = \\
&=&\left( -1\right) ^{n-1-7+1}\frac{(n-1+2)!\chi \left( {\cal M}%
_{2,7}\right) }{8!}= \\
&=&\left( -1\right) ^{n+1}(n+1)!\frac{\chi \left( {\cal M}_{2,7}\right) }{8!}%
\text{.}
\end{eqnarray*}

\noindent Moreover, for the case $n=7$, the fiber of the universal curve 
\begin{equation*}
{\cal M}_{2,7}\rightarrow {\cal M}_{2,6}\text{,}
\end{equation*}
over ${\cal M}_{2,6}\backslash U_{6,0}$ is a genus $2$ curve without $6$
points, of characteristic $-8$, and over $U_{6,0}$ is a genus $2$ curve
without the $6$ points fixed by the involution, modulo the involution
itself, hence a genus $0$ curve without $6$ points, of characteristic $-4$;
therefore, since $U_{6,0}\cong {\cal M}_{0,6}$ , we get 
\begin{eqnarray*}
\chi \left( {\cal M}_{2,7}\right) &=&-8\left( \chi \left( {\cal M}%
_{2,6}\right) -\chi \left( {\cal M}_{0,6}\right) \right) -4\chi \left( {\cal %
M}_{0,6}\right) + \\
&=&-8\chi \left( {\cal M}_{2,6}\right) -24\text{.}
\end{eqnarray*}

\noindent Since we have all the ingredients, we begin computing directly 
the cases $n=1,...,6$.

\begin{equation*}
\chi \left( {\cal M}_{2,0}\right) =\chi \left( \frac{{\cal M}_{0,6}}{S_6}%
\right) =1\text{.}
\end{equation*}
 
\begin{eqnarray*}
\chi \left( {\cal M}_{2,1}\right) =\sum_{j=0}^12^{-j}\chi \left( \frac{{\cal %
M}_{0,7-j}}{S_{6-j}}\right) +\frac 12\chi \left( \frac{{\cal M}_{0,6}}{S_5}%
\right)  \\
\ =\chi \left( \frac{{\cal M}_{0,7}}{S_6}\right) +\chi \left( \frac{{\cal M}%
_{0,6}}{S_5}\right) =2\text{.}
\end{eqnarray*}
  
\begin{eqnarray*}
\chi \left( {\cal M}_{2,2}\right) =2\chi \left( \frac{{\cal M}_{0,8}}{S_6}\right)
- \chi \left( \frac{{\cal M}_{0,5}}{S_3}\right) + \chi \left( \frac{{\cal M}%
_{0,7}}{S_6}\right) \\
+2\chi \left( \frac{{\cal M}_{0,7}}{S_5}\right) +\chi
\left( \frac{{\cal M}_{0,6}}{S_4}\right) 
 =2\text{.}
\end{eqnarray*}
  
\begin{eqnarray*}
&& \\
\chi \left( {\cal M}_{2,3}\right)=4\cdot \chi \left( \frac{{\cal M}_{0,9}}
{S_6}\right) 
+6\cdot \chi \left( 
\frac{{\cal M}_{0,8}}{S_6}\right) -3 \chi \left( \frac{{\cal M}_{0,5}}{S_3}\right) \\
+6\chi \left( \frac{{\cal M}_{0,8}}{S_5}%
\right) +3\chi \left( \frac{{\cal M}_{0,7}}{S_5}\right)  
+3\chi \left( \frac{{\cal M}_{0,7}}{S_4}\right) +\chi \left( \frac{%
{\cal M}_{0,6}}{S_3}\right) =0\text{.} \\
\end{eqnarray*}

\noindent With the same kind of computations, we get that  $\chi \left( {\cal M}_{2,4}\right)=-4$ , 
$ \chi \left( {\cal M}_{2,5}\right)=0 $, and
$\chi \left( {\cal M}_{2,6}\right)=-24$. 
 Finally, 
\begin{equation*}
\chi \left( {\cal M}_{2,7}\right) =-8\chi \left( {\cal M}_{2,6}\right)
-24=168\text{.}
\end{equation*}

From this we can conclude that, for $n\geq 7$, 
\begin{eqnarray*}
\chi \left( {\cal M}_{2,n}\right) &=&\left( -1\right) ^{n+1}(n+1)!\frac{\chi
\left( {\cal M}_{2,7}\right) }{8!} \\
&=&\left( -1\right) ^{n+1}\frac{(n+1)!}{240}\text{.}
\end{eqnarray*}

In the following table we list our results: 
\begin{equation*}
\begin{tabular}{|lllllllll|}
\hline
\multicolumn{1}{|l|}{$n$} & \multicolumn{1}{l|}{$0$} & \multicolumn{1}{l|}{$%
1 $} & \multicolumn{1}{l|}{$2$} & \multicolumn{1}{l|}{$3$} & 
\multicolumn{1}{l|}{$4$} & \multicolumn{1}{l|}{$5$} & \multicolumn{1}{l|}{$6$%
} & $\geq 7$ \\ \hline
\multicolumn{1}{|l|}{$\chi \left( {\cal M}_{2,n}\right) $} & 
\multicolumn{1}{l|}{$1$} & \multicolumn{1}{l|}{$2$} & \multicolumn{1}{l|}{$2$%
} & \multicolumn{1}{l|}{$0$} & \multicolumn{1}{l|}{$-4$} & 
\multicolumn{1}{l|}{$0$} & \multicolumn{1}{l|}{$-24$} & $\left( -1\right)
^{n+1}\frac{(n+1)!}{240}$ \\ \hline
\end{tabular}
\end{equation*}
\section{The Euler characteristic of ${\cal M}_{1,n}$}
\label{m1n}

A similar calculation can be done for the Euler characteristic of ${\cal M}%
_{1,n}$ (see \cite{AC}). In this case, $\tau $ denotes the hyperelliptic
involution around the last marked point, and the subvarieties $U_{j,r}$ have
the same definition, assuming that $j$ counts the $\tau $-fixed marked
points except the last one. In this case, for $n\geq 5$, the curves have no
automorphisms, and we get 
\begin{eqnarray*}
\chi \left( {\cal M}_{1,n}\right) &=&\prod_{h=5}^{n-1}-h\chi \left( {\cal M}%
_{1,5}\right) \\
\ &=&\left( -1\right) ^{n-1}\frac{(n-1)!\chi \left( {\cal M}_{1,5}\right) }{%
4!}\text{;}
\end{eqnarray*}
moreover, the fiber of the universal curve 
\begin{equation*}
{\cal M}_{1,5}\rightarrow {\cal M}_{1,4}\text{,}
\end{equation*}
over ${\cal M}_{1,4}\backslash U_{3,0}$ is a genus $1$ curve without $4$
points, of characteristic $-4$, and over $U_{3,0}$ a genus $0$ curve without 
$4$ points, of characteristic $-2$; we then get 
\begin{eqnarray*}
\chi \left( {\cal M}_{1,5}\right) &=&-4\left( \chi \left( {\cal M}%
_{1,4}\right) -\chi \left( {\cal M}_{0,4}\right) \right) -2\chi \left( {\cal %
M}_{0,4}\right) \\
\ &=&-4\chi \left( {\cal M}_{1,4}\right) -2\text{.}
\end{eqnarray*}

\noindent The formula we get for $n\leq 3$ is 
\begin{eqnarray*}
\chi \left( {\cal M}_{1,n+1}\right) &=&\sum_{j=0}^n\sum_{r=0}^{\left[ \frac{%
n-j}2\right] }a_{j,r}\chi \left( U_{j,r}\right) \\
\ &=&\sum_{j=0}^n\sum_{r=0}^{\left[ \frac{n-j}2\right] }\binom nj\frac{(n-j)!%
}{(n-j-2r)!r!}2^{n-j-2r-1}\chi \left( \frac{{\cal M}_{0,n-r+4-j}}{S_{3-j}}%
\right) + \\
&&\ +\frac 12\chi \left( \frac{{\cal M}_{0,4}}{S_{3-n}}\right) \text{ .}
\end{eqnarray*}

\noindent Since we know that ${\cal M}_{1,1}\cong {\Bbb C}$, then $\chi \left( {\cal M}%
_{1,1}\right) =1$ ; from the formula we calculate:

\begin{itemize}
\item  
\begin{eqnarray*}
\chi \left( {\cal M}_{1,2}\right) =\sum_{j=0}^12^{-j}\chi \left( \frac{{\cal %
M}_{0,5-j}}{S_{3-j}}\right) +\frac 12\chi \left( \frac{{\cal M}_{0,4}}{S_2}%
\right) =1\text{;}
\end{eqnarray*}

\item  
\begin{eqnarray*}
\chi \left( {\cal M}_{1,3}\right) =\sum_{j=0}^2\sum_{r=0}^{\left[ \frac{2-j}%
2\right] }\binom 2j\frac{(2-j)!}{(2-j-2r)!}2^{1-j-2r}\chi \left( \frac{{\cal %
M}_{0,6-r-j}}{S_{3-j}}\right) +\frac 12\chi \left( {\cal M}_{0,4}\right)  \\
\ =2\chi \left( \frac{{\cal M}_{0,6}}{S_3}\right) +3\chi \left( \frac{{\cal M%
}_{0,5}}{S_2}\right) +\chi \left( {\cal M}_{0,4}\right) =0\text{;}
\end{eqnarray*}

\item  
\begin{eqnarray*}
\chi \left( {\cal M}_{1,4}\right) =\sum_{j=0}^3\sum_{r=0}^{\left[ \frac{3-j}%
2\right] }\binom 3j\frac{(3-j)!}{(3-j-2r)!r!}2^{2-j-2r}\chi \left( \frac{%
{\cal M}_{0,7-r-j}}{S_{3-j}}\right) +\frac 12\chi \left( {\cal M}%
_{0,4}\right)  \\
\ =\sum_{j=0}^3\sum_{r=0}^{\left[ \frac{3-j}2\right] }\binom 3j\frac{%
(3-j)!2^{2-j-2r}\left( -1\right) ^{r+j}\left( 4-r-j\right) !}{%
(3-j-2r)!r!\left( 3-j\right) !}-\frac 12=0\text{.}
\end{eqnarray*}
\end{itemize}

Finally, $\chi \left( {\cal M}_{1,5}\right) =-2$ and, for $n\geq 5$, $\chi
\left( {\cal M}_{1,n}\right) =\left( -1\right) ^n\frac{(n-1)!}{12}$; in the
following table we summarize the results: 
\begin{equation*}
\begin{tabular}{|llllll|}
\hline
\multicolumn{1}{|l|}{$n$} & \multicolumn{1}{l|}{$1$} & \multicolumn{1}{l|}{$%
2 $} & \multicolumn{1}{l|}{$3$} & \multicolumn{1}{l|}{$4$} & $\geq 5$ \\ 
\hline
\multicolumn{1}{|l|}{$\chi \left( {\cal M}_{1,n}\right) $} & 
\multicolumn{1}{l|}{$1$} & \multicolumn{1}{l|}{$1$} & \multicolumn{1}{l|}{$0$%
} & \multicolumn{1}{l|}{$0$} & $\left( -1\right) ^n\frac{(n-1)!}{12}$ \\ 
\hline
\end{tabular}
\end{equation*}
\begin{section}{Generating functions and graphs}

\label{graphs}

In this section we briefly describe our combinatorial strategy in
computing the Euler characteristics of \( {\overline {\cal M}}_{1,n} \)
and \({\overline {\cal
M}}_{2,n} \).

As mentioned in the Introduction, the quasi-projective subvarieties \( X_i
\) of \({\overline {\cal M}}_{g,n} \) are in correspondence with a collection of
genus \( g \) graphs with \( n \) leaves. Considering a graph \( \Gamma
\) of this collection, its contribution to the Euler characteristic is
provided by the product of the Euler characteristics of the moduli spaces
associated to its  vertices.

Now it turns out that all the genus \( g \) graphs (with any number of
leaves) can be obtained by attaching trees to some loops or to some
vertex representing a curve of genus greater than 1.
Therefore, if we let \( D \) to be a generating function  which counts
the contributions of trees to the Euler characteristic, the computation
of the series \( K_g \) reduces to the sum of the contributions of all
the possible combinations of loops and  vertices with genus greater than
1, multiplied by a suitable power of \( D \).

Let us then give the following suitable definition of \( D \):
\begin{definition}
\[ D(t)=t+\sum_{n=2}^{\infty}\sum_{\Gamma}\chi (\Gamma)\frac{t^n}{n!} \]
Here the second sum ranges over the collection of ``admissible trees'',
which are all the oriented rooted trees on \( n \) numbered
leaves with the following two further properties
\begin{enumerate}
        \item  There is an unlabelled half-edge going into the root.

        \item  The tree is stable, that is to say, for every  vertex of
        the tree the number of outgoing edges (including the leaves),
        plus 1, is greater than or equal to 3.
\end{enumerate}
Furthermore, \( \chi (\Gamma) \) is the Euler characteristic of the stratum
of \( {\cal M}_{0,n} \) which corresponds to \( \Gamma \).
\end{definition}

As an example of the way in which we use \( D \), the contribution to \(
K_2 \) of all the graphs of the following kind

\begin{figure}[h]
\begin{center} 
\mbox{\epsfig{file=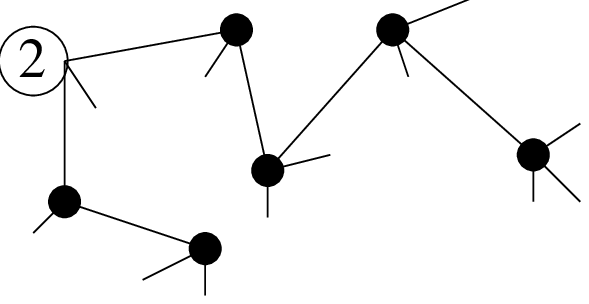,width=2cm,height=5cm,angle=270}}
\label{fun1}
\end{center}
\end{figure}

can be computed by means of the series
\begin{displaymath}
        \sum_{n\geq 1}\chi ({\cal M}_{2,n})\frac{D^n}{n!}
\end{displaymath}

Let us now focus on the series \( D \): we can easily find a recursive
relation for it  noticing that, if we cut an edge which stems from the
root of a genus 0 admissible tree, the cut part is again a genus 0
 admissible tree. 
 
\noindent Therefore we can write the following recursive relation for $D$ (which
was taken as a definition in the Introduction):
\[ D=t+\sum_{n=2}^{\infty}\chi  ({\cal M}_{0,n+1})\frac{D^n}{n!} \text{.} \]
Substituting the values for \(\chi  ({\cal M}_{0,n+1})  \) we obtain
\[ D=t+\sum_{n=2}^{\infty}(-1)^{n-1}(n-2)!\frac{D^n}{n!} \text{,}\]
which, after differentiating with respect to \( t \) gives
\[ D'=1+D'\sum_{n=1}^{\infty}(-1)^{n}\frac{D^n}{n} \]

\[ D'(1-log(1+D))=1 \text{.} \]
This is a differential equation that allows us to compute  recursively
all the coefficients of \( D \text{:} \)
\[ D(t)= t + \frac{t^2}2 + \frac{t^3}3 + \frac{7t^4}{24} + \frac{17t^5}{60} + 
    \frac{71t^6}{240} + \frac{163t^7}{504} +o\left( t^8 \right)  \text{.}       \]

We notice that, in computing  Euler characteristics,
we will often come across the series \(E = log(1+D) \), which can be
equivalently
written in the following way
\begin{displaymath}
       E=\sum_{n=1}^{\infty}\chi  ({\cal M}_{0,n+2})\frac{D^n}{n!}
\end{displaymath}
and summarizes the contribution provided by trees that
stem from a vertex of a polygon.

Let us now pass to motivate the introduction of two different
operations on generating series, namely the derivative with respect to \(
D \) and the derivative with respect to \( t \).
Let us consider a graph $ \Delta $ of the following kind

\newpage

\begin{figure}[h]
\begin{center} 
\mbox{\epsfig{file=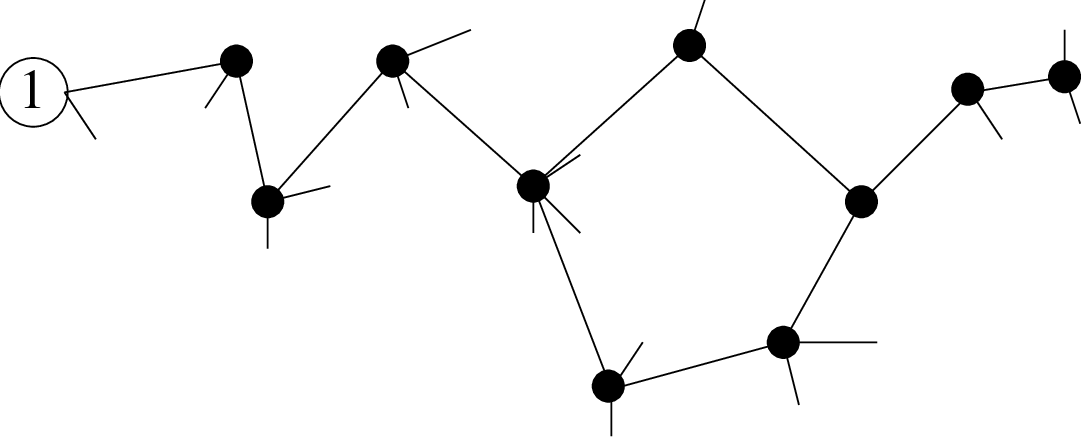,width=2cm,height=5cm,angle=270}}
\label{fun2}
\end{center}
\end{figure}

\noindent or more generally a graph which is made by two components of genus 1
attached by a genus 0 path.

When computing the contribution of such graphs to the series \( K_2 \),
we can imagine to cut the picture in the following way

\begin{figure}[h]
\begin{center} 
\mbox{\epsfig{file=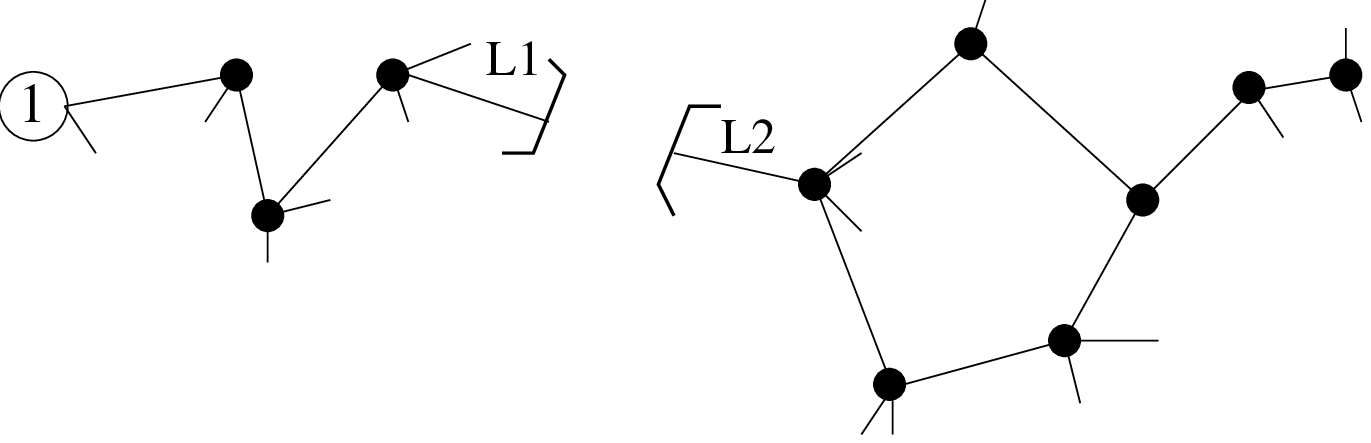,width=2cm,height=6cm,angle=270}}
\label{fun3}
\end{center}
\end{figure}

What remains on the left is a  graph \( \Delta_1 \) of genus 1 with the
artificial leaf \( L_1 \), while on the right there is another graph \(
\Delta_2 \) of
genus 1 with the artificial leaf \( L_2 \) instead of a tree.

A first approximation to the contribution of such graphs is provided by
\(\frac{1}{2}\frac{\partial K_1}{\partial t} \frac{\partial K_1}{\partial
D} \). In fact
\begin{enumerate}
        \item  The factor \( \frac{\partial K_1}{\partial t} \) takes into
        account the contribution of the left part of the graph.  The
        derivative with respect to \( t \) cancels
        the mistake due to the presence in \( \Delta_1 \)  of the leaf \( L_1 \)
        which is not a leaf
        of \( \Delta  \).

        \item The factor \( \frac{\partial K_1}{\partial D } \) takes into
        account the contribution of the right part of the graph. The
derivative with
respect to \( D \)
        is the translation in
        terms of generating series of the presence, in  \( \Delta_2 \),
        of the artificial leaf \( L_2 \) instead of a tree.

        \item The coefficient \( \frac{1}{2} \) is needed since in general
there are two
        possible ways to cut the graph \( \Delta \).
\end{enumerate}

There are some exceptional cases; for instance one is provided by
the following
configuration, which should be carefully looked at:

\begin{figure}[h]
\begin{center} 
\mbox{\epsfig{file=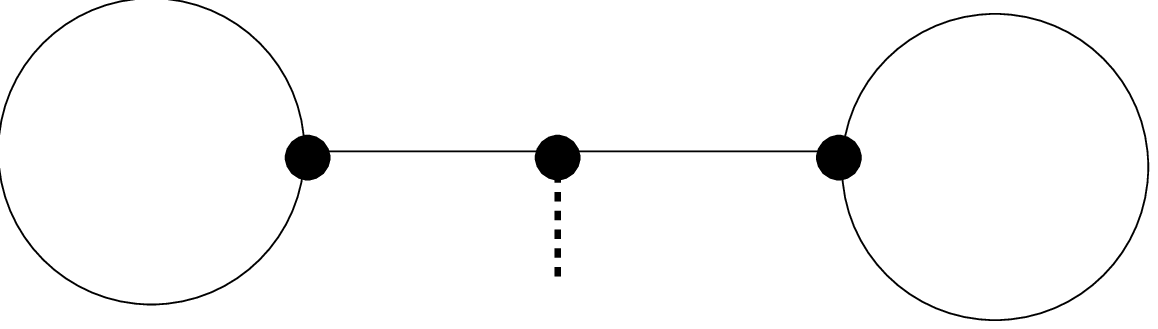,width=2cm,height=6cm,angle=270}}
\label{fun4}
\end{center}
\end{figure}

Here and from now on, when we draw a dotted line outgoing from a vertex 
we mean that any number of trees can stem from the vertex itself.
 
In this  case, the  contribution to \( K_2 \) turns out to be equal to
\begin{displaymath}
        \sum_{n \geq 1}\chi  ({\cal M}_{0,n}/S_2)\frac{D^n}{n!} 
\end{displaymath}
instead of
\begin{displaymath}
        \frac{1}{2}\sum_{n \geq 1}\chi  ({\cal M}_{0,n})\frac{D^n}{n!}
\end{displaymath}
which was implicit in the expression \(\frac{1}{2}\frac{\partial
K_1}{\partial t} \frac{\partial K_1}{\partial
D} \).

\end{section}

\section{The generating function for genus $1$}
\label{k1}

In this section we obtain with our elementary and direct methods the
generating function for the genus one case, which was already calculated by
Getzler in \cite{G1}.

There are two different contributions which should be considered: the first
one comes from graphs of this kind

\begin{figure}[h]
\begin{center} 
\mbox{\epsfig{file=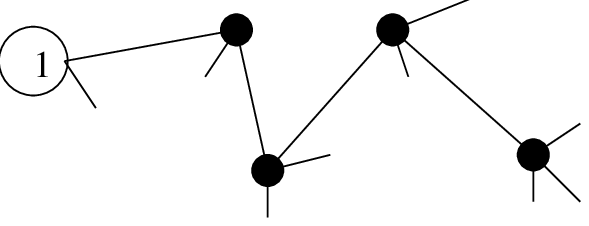,width=1.3cm,height=4cm,angle=270}}
\label{gen1}
\end{center}
\end{figure}

\noindent and is exactly 
\begin{eqnarray*}
\sum_{n\geq 1}\chi \left( {\cal M}_{1,n}\right) \frac{D^n}{n!} &=&D+\frac{D^2%
}2+\frac 1{12}\left( \sum_{n\geq 1}\left( -1\right) ^n\frac{D^n}n+D-\frac{D^2%
}2+\frac{D^3}3-\frac{D^4}4\right) \\
&=&\frac{13}{12}D+\frac{11}{24}D^2+\frac{D^3}{36}-\frac{D^4}{48}-\frac
1{12}\log \left( 1+D\right) \\
&=&\frac{13}{12}D+\frac{11}{24}D^2+\frac{D^3}{36}-\frac{D^4}{48}-\frac E{12}%
\text{;}
\end{eqnarray*}
the second contribution comes from the graphs containing a loop, and is 
\begin{equation*}
\frac 12\sum_{l\geq 3}\frac{E^l}l+\sum_{n\geq 1}\chi \left( \frac{{\cal M}%
_{0,n+2}}{S_2}\right) \frac{D^n}{n!}+\frac 12\sum_{n,m\geq 1}\chi \left( 
\frac{{\cal M}_{0,n+2}\times {\cal M}_{0,m+2}}{S_2}\right) \frac{D^{n+m}}{%
n!m!}\text{,}
\end{equation*}
where the action of $S_2$ exchanges the last two markings, simultaneously in
the second case.

It can be written as 
\begin{eqnarray*}
&&\ \frac 12\left( \log \left( 1-E\right) -E-\frac{E^2}2\right) +\frac
12\left( E-D+\frac{D^2}2\right) +\frac 14\left( E^2-D^2+D^3-\frac{D^4}%
4\right) \\
&&\ +\sum_{1\leq n\leq 2}\chi \left( \frac{{\cal M}_{0,n+2}}{S_2}\right) 
\frac{D^n}{n!}+\frac 12\sum_{1\leq n,m\leq 2}\chi \left( \frac{{\cal M}%
_{0,n+2}\times {\cal M}_{0,m+2}}{S_2}\right) \frac{D^{n+m}}{n!m!} \\
\ &=&\frac 12\left( -\log \left( 1-E\right) -E-\frac{E^2}2\right) +\frac
12\left( E-D+\frac{D^2}2\right) + \\
&&\ +\frac 14\left( E^2-D^2+D^3-\frac{D^4}4\right) +D+\frac{D^2}2+\frac{D^4}8
\\
\ &=&-\frac 12\log \left( 1-E\right) +\frac D2+\frac{D^2}2+\frac{D^3}4+\frac{%
D^4}{16}\text{.}
\end{eqnarray*}
Summing up we get the generating function for genus $1$. 
\begin{eqnarray*}
K_1 &=&\frac{13}{12}D+\frac{11}{24}D^2+\frac{D^3}{36}-\frac{D^4}{48}-\frac
E{12}-\frac 12\log \left( 1-E\right) +\frac D2+\frac{D^2}2+\frac{D^3}4+\frac{%
D^4}{16} \\
\ &=&\frac{19}{12}D+\frac{23}{24}D^2+\frac 5{18}D^3+\frac{D^4}{24}-\frac
E{12}-\frac 12\log \left( 1-E\right) \text{.}
\end{eqnarray*}

\section{Graphs and components} 
\label{k2}

Let us consider the set ${\cal F}$ of graphs representing boundary
components of the moduli space of genus $2$ pointed curves, including the
one representing the open part. We define the {\bf graph-type} of each graph
as the graph obtained by deleting the markings, contracting all trees, and
smoothing each vertex of valence $2$.

There are seven different graph-types, and we denote by ${\cal F}_1,...,%
{\cal F}_7$ the collections of graphs of each type. We count separately the
contribution these sets give to the generating function for the
characteristic of $\overline{{\cal M}}_{2,n}$.

\begin{figure}[h]
\begin{center} 
\mbox{\epsfig{file=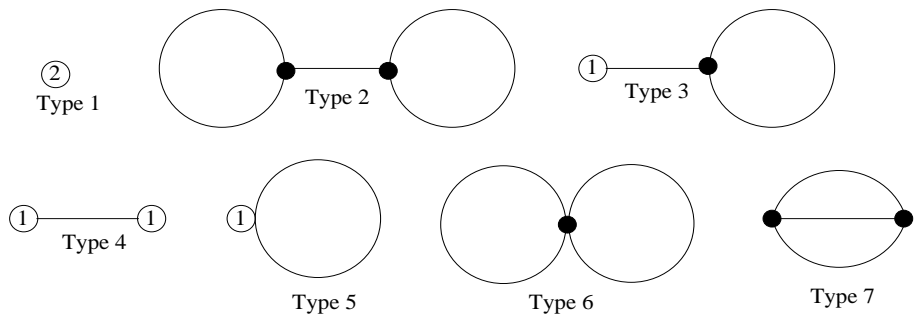,width=4cm,height=12cm,angle=270}}
\caption[]{Graph-types}
\label{types}
\end{center}
\end{figure}

The group of automorphisms $G_i$ of the graph-type acts on ${\cal F}_i$, and
generically the contribution should have coefficient $\frac 1{|G_i|}$ ; in
this general formula we should ``correct'' the contribution of the
components corresponding to graphs with non trivial stabilizer. In this
case, the stabilizer of the graph acts generically non-trivially on the
boundary component, since it acts on the added markings of the irreducible
components of the curve. If there are enough markings on this component ,
the action is free, and nothing changes in the coefficient $\frac 1{|G_i|}$;
but, for a low number of markings, we should analyze the contribution this
component gives independently.

As an example, graphs of this kind

\begin{figure}[h]
\begin{center} 
\mbox{\epsfig{file=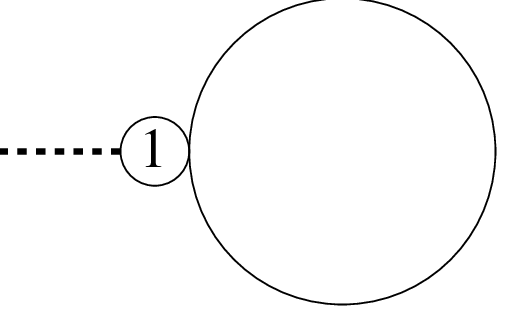,width=2cm,height=3cm,angle=270}}
\label{graph1}
\end{center}
\end{figure}

\noindent (which are of type $5$) give the contribution $\frac 12\sum_{n\geq
0}\chi \left( {\cal M}_{1,n+2}\right) \frac D{n!}^n$ in the generic formula,
since the graph-type has automorphism group of order $2$. In fact, what we
should really put in the formula is $\sum_{n\geq 0}\chi \left( \frac{{\cal M}%
_{1,n+2}}{S_2}\right) \frac D{n!}^n$, where $S_2$ permutes the two added
markings. The two formulas coincides for $n\geq 5$, when the action becomes
free, but the initial terms should be corrected.

Now let us analyze separately each contribution.

\subsection{Graphs of type 1}

Since there are no graph-type automorphisms, the contribution is $%
\sum_{n\geq 0}\chi \left( {\cal M}_{2,n}\right) \frac{D^n}{n!}$.

\noindent From our formulas, we get 
\begin{eqnarray*}
\sum_{n\geq 0}\chi \left( {\cal M}_{2,n}\right) \frac{D^n}{n!} &=&1+2D+D^2
-\frac{D^4}{6}-\frac{D^6}{30}+\frac 1{240}\sum_{n\geq
7}\left( -1\right) ^{n+1}(n+1)D^n= \\
&=&-\frac 1{240\left( 1+D\right) ^2}+\frac{241}{240}+\frac{239D}{120}+\frac{%
81D^2}{80}-\frac{D^3}{60}-\frac{7D^4}{48}-\frac{D^5}{40}-\frac{D^6}{240}%
\text{.}
\end{eqnarray*}
\bigskip

\subsection{Graphs of type 2,3,4}

\ These graph-types could be considered together, as they all could be seen
in exactly two different ways as the union of a genus $1$ graph with one cut
leaf, and a genus $1$ graph with one cut tree.

We recall that the generating function for the genus $1$ case is 
\begin{equation*}
K_1\left( D\right) =\frac{19D}{12}+\frac{23D^2}{24}+\frac{5D^3}{18}+\frac{D^4%
}{24}-\frac{\log \left( 1+D\right) }{12}-\frac{\log \left( 1-\log \left(
1+D\right) \right) }2\text{;}
\end{equation*}
the generic contribution is: 
\begin{eqnarray*}
&&\ \frac 12\left( \frac{\partial K_1\left( D\right) }{\partial D}\right)
^2D^{\prime } \\
\  &=&\frac 12\left( \frac{19}{12}+\frac{23D}{12}+\frac{5D^2}6+\frac{D^3}%
6-\frac 1{12\left( 1+D\right) }+\frac 1{2\left( 1-E\right) \left( 1+D\right)
}\right) ^2\frac 1{1-E}
\end{eqnarray*}
\begin{eqnarray*}
\  &=&\frac{\left( 361+874D+909D^2+536D^3+192D^4+40D^5+4D^6\right) }{%
288\left( 1-E\right) } \\
&&\ -\frac{19+23D+10D^2+2D^3}{144\left( 1-E\right) \left( 1+D\right) }+\frac{%
19+23D+10D^2+2D^3}{24\left( 1-E\right) ^2\left( 1+D\right) }-\frac
1{24\left( 1-E\right) ^2\left( 1+D\right) ^2} \\
&&\ +\frac 1{288\left( 1-E\right) \left( 1+D\right) ^2}+\frac 1{8\left(
1-E\right) ^3\left( 1+D\right) ^2}\text{;}
\end{eqnarray*}
moreover, the group exchanging the two sides of the graph type fixes the
following graphs of type $2$ and $4$:

\begin{figure}[h]
\begin{center} 
\mbox{\epsfig{file=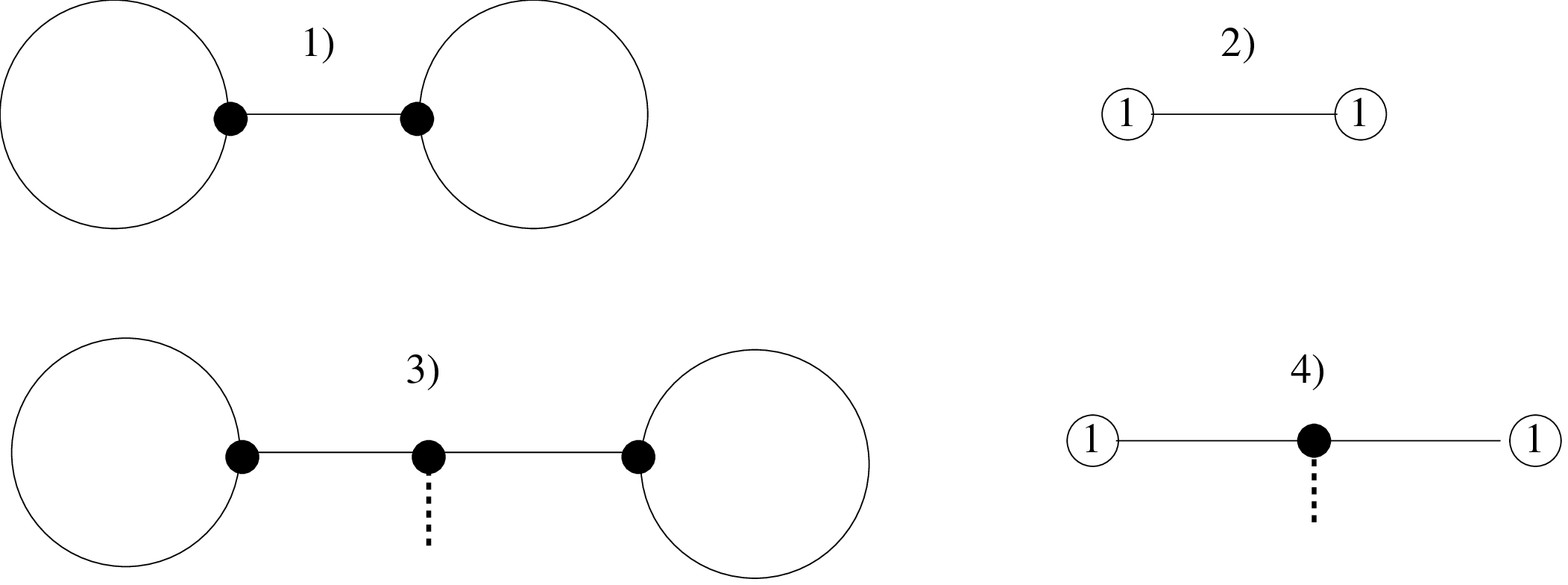,width=4cm,height=10cm,angle=270}}
\label{graph2}
\end{center}
\end{figure}

\noindent in each of the first two cases, the contribution is $1$ instead of 
$\frac 12$, while in the third one we should replace $\frac 12\left( D-\frac{D^2%
}2\right) $ with 
\begin{equation*}
\sum_{1\leq n\leq 2}\chi \left( \frac{{\cal M}_{0,n+2}}{S_2}\right) \frac{D^n%
}{n!}=D\text{.}
\end{equation*}

The fourth graph requires a little more work: in fact we have to analyze the action
of $S_2$ on  
$${\cal M}_{0,n+2} \times {\cal M}_{1,1} \times {\cal M}_{1,1} \text{,}$$ which is free for
$n>2$.
For the case $n=1$, the action is trivial, whereas for the case $n=2$, the quotient map 
is generically $2:1$, and it is ramified on 
$$\left[ {\Bbb P}^{1}; \infty, 0, 1,-1 \right] \times \Delta \text{,}$$ where $\Delta$ is 
the diagonal in ${\cal M}_{1,1} \times {\cal M}_{1,1}$.  Thus we get 
$$\chi \left( {\cal M}_{0,n+2} \times {\cal M}_{1,1} \times {\cal M}_{1,1} / S_2 \right) =0 \text{.}$$
Once more we should replace $\frac 12\left( D-\frac{D^2}2\right) $ with $D$.

Finally we have: 
\begin{eqnarray*}
&&\ 1+D+\frac{D^2}2+\frac{\left(
361+874D+909D^2+536D^3+192D^4+40D^5+4D^6\right) }{288\left( 1-E\right) } \\
&&\ \ -\frac{19+23D+10D^2+2D^3}{144\left( 1-E\right) \left( 1+D\right) }+%
\frac{19+23D+10D^2+2D^3}{24\left( 1-E\right) ^2\left( 1+D\right) }-\frac
1{24\left( 1-E\right) ^2\left( 1+D\right) ^2} \\
&&\ \ +\frac 1{288\left( 1-E\right) \left( 1+D\right) ^2}+\frac 1{8\left(
1-E\right) ^3\left( 1+D\right) ^2}\text{.}
\end{eqnarray*}
\bigskip

\subsection{Graphs of type 5}

\ $G_5\cong S_2$ acts on the graph-type. The formula to be corrected is 
\begin{eqnarray*}
&&\frac 12\left( \sum_{n\geq 0}\chi \left( {\cal M}_{1,n+2}\right) \frac{D^n%
}{n!}\right) D^{\prime } \\
&=&\frac 1{2\left( 1-E\right) }\left( \sum_{n\geq 3}\left( -1\right) ^n\frac{%
\left( n+1\right) !}{12}\frac{D^n}{n!}+1\right) \\
\ &=&\frac 1{2\left( 1-E\right) }\left( \frac 1{12}\left( \sum_{n\geq
0}\left( -1\right) ^n\left( n+1\right) D^n\right) -\frac 1{12}+\frac D6-%
\frac{D^2}4+1\right) \\
&=&\frac 1{2\left( 1-E\right) }\left( \frac 1{12\left( 1+D\right) ^2}+\frac{%
11}{12}+\frac D6-\frac{D^2}4\right)
\end{eqnarray*}
the graphs stabilized by $G_5$ are

\begin{figure}[h]
\begin{center} 
\mbox{\epsfig{file=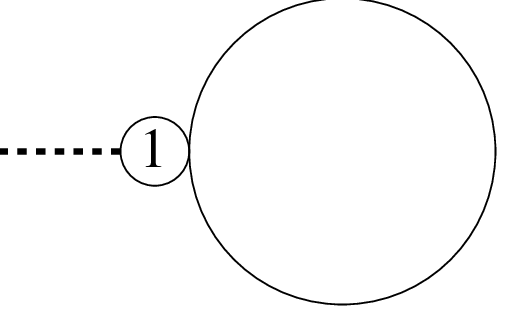,width=2cm,height=3cm,angle=270}}
\label{graph3}
\end{center}
\end{figure}

\noindent and

\begin{figure}[h]
\begin{center} 
\mbox{\epsfig{file=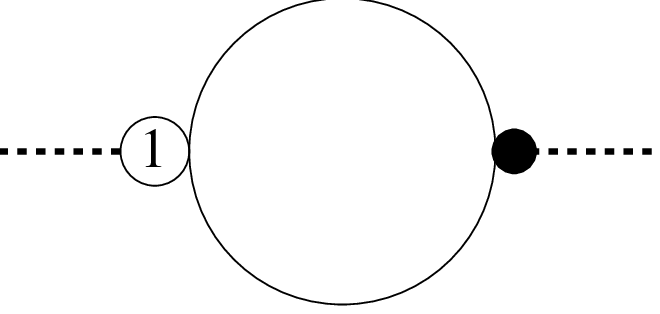,width=1.5cm,height=3cm,angle=270}}
\label{graph4}
\end{center}
\end{figure}

\noindent which, taking into account only the cases where the elliptic curve
has non trivial automorphisms, give contribution 
\begin{eqnarray*}
&&\ \ \ \ \sum_{0\leq n\leq 4}\chi \left( \frac{{\cal M}_{1,n+2}}{S_2}%
\right) \frac{D^n}{n!} \\
\ &=&\chi \left( \frac{{\cal M}_{1,2}}{S_2}\right) +\chi \left( \frac{{\cal M%
}_{1,3}}{S_2}\right) D+\chi \left( \frac{{\cal M}_{1,4}}{S_2}\right) \frac{%
D^2}2 \\
&&\ \ +\chi \left( \frac{{\cal M}_{1,5}}{S_2}\right) \frac{D^3}6+\chi \left( 
\frac{{\cal M}_{1,6}}{S_2}\right) \frac{D^4}{24} \\
\ &=&1+D+\frac{D^2}2 +\frac{D^4}4 
\end{eqnarray*}
instead of 
\begin{equation*}
\ \ \frac 12\sum_{0\leq n\leq 4}\chi \left( {\cal M}_{1,n+2}\right) \frac{D^n%
}{n!}=\frac 12\left( 1-\frac{D^3}3+\frac{5D^4}{12}\right) \text{,}
\end{equation*}
and 
\begin{equation*}
\ \sum\begin{Sb} 0\leq n\leq 4  \\ 1\leq m\leq 2  \end{Sb}  \chi \left( 
\frac{{\cal M}_{1,n+2}\times {\cal M}_{0,m+2}}{S_2}\right) \frac{D^{n+m}}{%
n!m!}=D+D^2+D^3+\frac{D^4}{4}+\frac{D^5}6 -\frac{D^6}{12}
\end{equation*}
instead of 
\begin{eqnarray*}
&&\ \ \ \ \frac 12\sum_{0\leq n\leq 4}\chi \left( {\cal M}_{1,n+2}\right) 
\frac{D^n}{n!}\left( D-\frac{D^2}2\right) \\
\ &=&\frac 12\left( 1-\frac{D^3}3+\frac{5D^4}{12}\right) \left( D-\frac{D^2}%
2\right) \\
\ &=&\frac D2-\frac{D^2}4-\frac{D^4}6+\frac{7D^5}{24}-\frac{5D^6}{48}\text{.}
\end{eqnarray*}

Therefore the contribution of graphs of type $5$ is 

\begin{eqnarray*}
&=&\frac 1{24\left( 1-E\right) \left( 1+D\right) ^2}+\frac{11+2D-3D^2}{%
24\left( 1-E\right) } \\
&&+\frac 12+\frac 32D+\frac 74D^2+\frac 76D^3+\frac{11}{24}D^4-\frac{1}{8}D^5+%
\frac{1}{48}D^6\text{.}
\end{eqnarray*}
\bigskip

\subsection{Graphs of type 6}

\ $G_6=H_7\times L_7\times R_7\cong S_2\times S_2\times S_2$ acts on the
graph type: $H_7$ exchanges the two loops, $L_7$ and $R_7$ reverse the
orientation on the left and on the right loop; the starting formula is 
\begin{equation*}
\frac 18\sum_{n\geq 0}\chi \left( {\cal M}_{0,n+4}\right) \frac{D^n}{n!}%
\left( D^{\prime }\right) ^2=-\frac 1{8\left( 1+D\right) ^2\left( 1-E\right)
^2}\text{.}
\end{equation*}

\begin{figure}[h]
\begin{center} 
\mbox{\epsfig{file=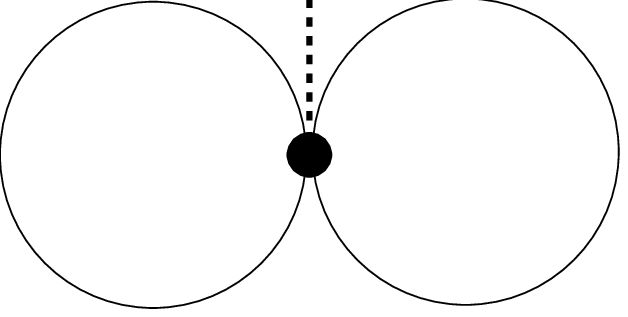,width=1.3cm,height=3cm,angle=270}}
\label{graph5}
\end{center}
\end{figure}
$G_6$ stabilizes the graph, and the contribution is 
\begin{eqnarray*}
&&\sum_{0\leq n\leq 3}\chi \left( \frac{{\cal M}_{0,n+4}}{D_4}\right) \frac{%
D^n}{n!} \\
&=&\chi \left( \frac{{\cal M}_{0,4}}{D_4}\right) +\chi \left( \frac{{\cal M}%
_{0,5}}{D_4}\right) D+\chi \left( \frac{{\cal M}_{0,6}}{D_4}\right) \frac{D^2%
}2 \\
&=&-\frac{D^2}2\text{,}
\end{eqnarray*}
which should replace $\frac 18\left( -1+2D-3D^2\right) $ in the formula.

\begin{figure}[h]
\begin{center} 
\mbox{\epsfig{file=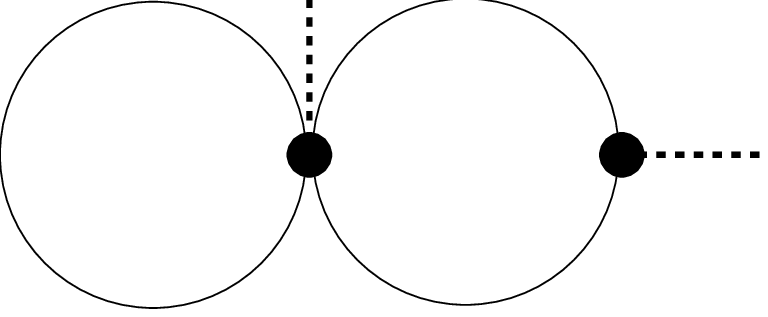,width=1.3cm,height=3cm,angle=270}}
\label{graph6}
\end{center}
\end{figure}

The stabilizer is $L_7\times R_7$; we put 
\begin{eqnarray*}
&&\frac 12 \chi \left( \frac{{\cal M}_{0,4}}{R_7}%
\right) \left( E-D+\frac{D^2}2\right) +D\sum_{0\leq n\leq
2}\chi \left( \frac{{\cal M}_{0,n+4}\times {\cal M}_{0,3}}{L_7\times R_7}%
\right) \frac{D^n}{n!} \\
&&+\frac{D^2}2\sum_{0\leq n\leq 2}\chi \left( \frac{{\cal M}_{0,n+4}\times 
{\cal M}_{0,4}}{L_7\times R_7}\right) \frac{D^n}{n!} \\
&=&-\frac 32{D^3}+\frac{D^4}4
\end{eqnarray*}
instead of 
\begin{eqnarray*}
&&\left( \frac D4 - \frac {D^2}8 \right) 
\left( -1+2D-3D^2\right) -\frac 14 \left( E-D+\frac{D^2}2 \right)  \\
&=&-\frac E4 +\frac{D^2}2 - D^3 +\frac 38{D^4} \text{.}
\end{eqnarray*}

\begin{figure}[h]
\begin{center} 
\mbox{\epsfig{file=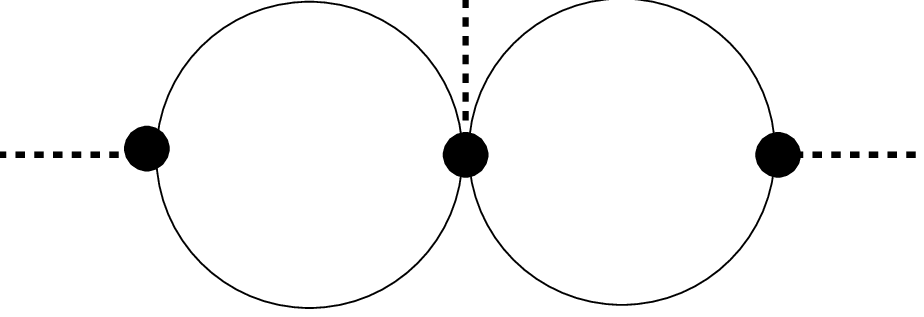,width=1.3cm,height=4cm,angle=270}}
\label{graph7}
\end{center}
\end{figure}

This set of graphs gives 
\begin{eqnarray*}
&&\ \ \ \frac 12\sum\begin{Sb} 0\leq n\leq 2  \\ 1\leq p,q\leq 2  \end{Sb} 
\chi \left( \frac{{\cal M}_{0,n+4}\times {\cal M}_{0,p+2}\times {\cal M}%
_{0,q+2}}{L_7\times R_7}\right) \frac{D^{n+p+q}}{n!p!q!} \\
\ &=&\frac 12\left( -\frac{D^4}4-D^4-D^4+\frac{D^5}2-\frac{D^6}4\right) \\
\ &=&-\frac 98D^4+\frac{D^5}4-\frac{D^6}8
\end{eqnarray*}
replacing 
\begin{eqnarray*}
&&\ \ \frac 18\left( -1+2D-3D^2\right) \left( D-\frac{D^2}2\right) ^2 \\
\ &=&-\frac{D^2}8+\frac{3D^3}8-\frac{21}{32}D^4+\frac 7{16}D^5-\frac 3{32}D^6
\end{eqnarray*}
and 
\begin{eqnarray*}
&&\frac 12\sum_{1\leq p\leq 2}\chi \left( \frac{{\cal M}_{0,4}\times {\cal M}%
_{0,p+2}}{S_2}\right) \frac{D^p}{p!}\left( E-D+\frac{D^2}2\right) \\
&=&\frac{ED^2}4-\frac{D^3}4+\frac{D^4}8
\end{eqnarray*}
replacing 
\begin{eqnarray*}
&&-\frac 14\left( E-D+\frac{D^2}2\right) \left( D-\frac{D^2}2\right) \\
&=&E\left( -\frac D4+\frac{D^2}8\right) +\frac{D^2}4-\frac{D^3}4+\frac{D^4}{%
16}\text{.}
\end{eqnarray*}

\begin{figure}[h]
\begin{center} 
\mbox{\epsfig{file=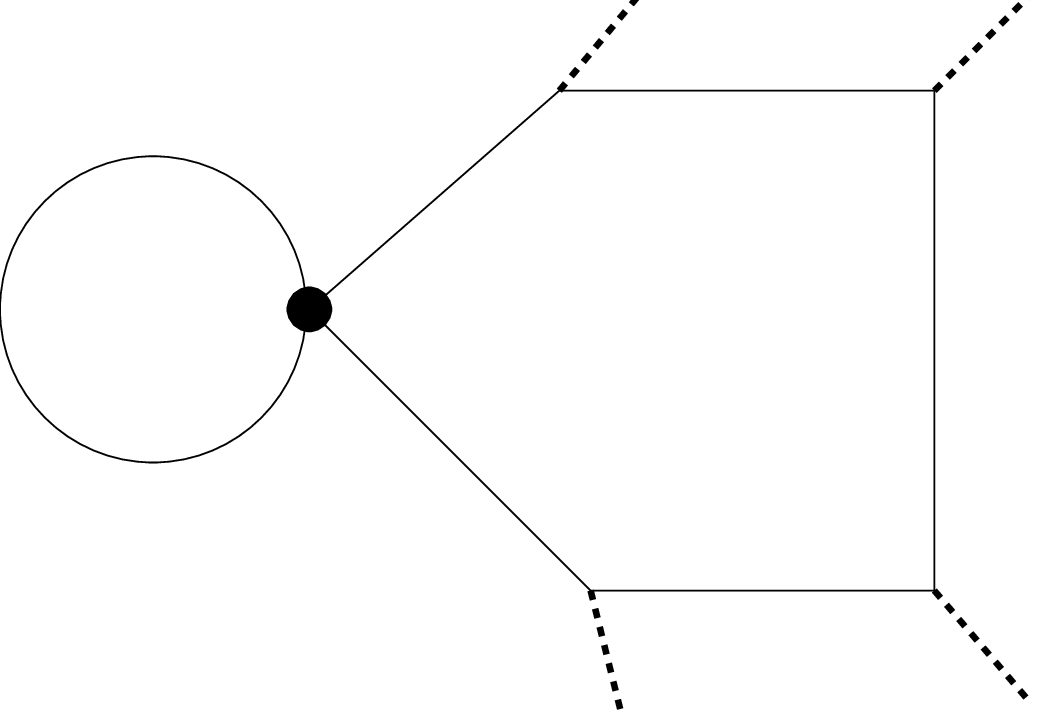,width=3cm,height=5cm,angle=270}}
\label{graph8}
\end{center}
\end{figure}

The contribution of these graphs (here and from now in the pictures a pentagon stands for
any polygon with more than three edges) is 
\begin{equation*}
\frac 12\chi \left( \frac{{\cal M}_{0,4}}{S_2}\right) \frac{E^2}{1-E}=0
\end{equation*}
where in the formula we have $-\frac 14\frac{E^2}{1-E}$,
and finally we have

\begin{figure}[h]
\begin{center} 
\mbox{\epsfig{file=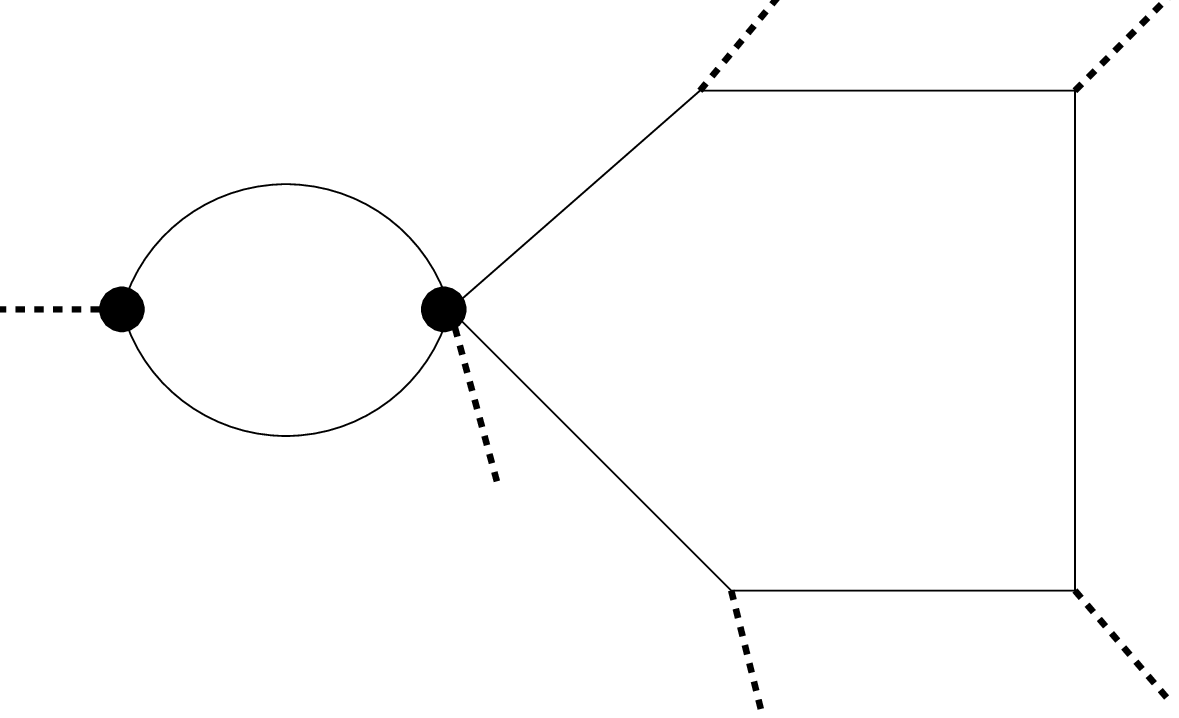,width=3cm,height=5cm,angle=270}}
\label{graph9}
\end{center}
\end{figure}

which give contribution 
\begin{eqnarray*}
&&\frac{E^2}{2\left( 1-E\right) }\left( \chi \left( \frac{{\cal M}%
_{0,4}\times {\cal M}_{0,3}}{S_2}\right) D+\chi \left( \frac{{\cal M}%
_{0,4}\times {\cal M}_{0,4}}{S_2}\right) \frac{D^2}2\right) \\
&=&\frac 14\frac{D^2E^2}{1-E}
\end{eqnarray*}
and not $\frac 14\frac{E^2}{1-E}\left( -D+\frac{D^2}2\right) $.

From graph type $6$:

\begin{eqnarray*}
&&\ -\frac 1{8\left( 1+D\right) ^2\left( 1-E\right) ^2}-\frac{D^2}2+\frac
18-\frac D4+\frac 38D^2 \\
&&\ -\frac32 D^3 +\frac{D^4}{4}+\frac{E}{4}-\frac{D^2}{2}+D^3-\frac38 D^4  \\
&&\ -\frac 98D^4+\frac{D^5}4-\frac{D^6}8+\frac{D^2}8-\frac{3D^3}8+\frac{21}{%
32}D^4-\frac 7{16}D^5+\frac 3{32}D^6 \\
&&\ +\frac{ED^2}4-\frac{D^3}4+\frac{D^4}8-E\left( -\frac D4+\frac{D^2}%
8\right) -\frac{D^2}4+\frac{D^3}4-\frac{D^4}{16} \\
&&\ +\frac 14\frac{E^2}{1-E}+\frac 14\frac{E^2D^2}{1-E}+\frac 14\frac{E^2D}{%
1-E}-\frac 18\frac{E^2D^2}{1-E} \\
\  &=&-\frac 1{8\left( 1+D\right) ^2\left( 1-E\right) ^2}+\frac{E^2}{4\left(
1-E\right) }\left( \frac{D^2}2+D+1\right)  \\
&&\ +E\left( \frac{D^2}8+\frac D4+\frac 14\right) +\frac 18-\frac D4-\frac{%
D^2}2-\frac 78D^3-\frac{11}{32}D^4-\frac 3{16}D^5-\frac{D^6}{32} \text{.}
\end{eqnarray*}
\bigskip

\subsection{Graphs of type 7}

The group acting is $G_7=H_7\times K_7$, where $H_7\cong S_2$ interchanges
the two vertices, and $K_7\cong S_3$ permutes the three edges.

The generic contribution is 
\begin{equation*}
\frac 1{12}\left( \sum_{n\geq 0}\chi \left( {\cal M}_{0,n+3}\right) \frac{D^n%
}{n!}\right) ^2\left( D^{\prime }\right) ^3=\frac 1{12}\frac 1{\left(
1+D\right) ^2}\frac 1{\left( 1-E\right) ^3}\text{;}
\end{equation*}
we need to correct the following contributions :

\begin{figure}[h]
\begin{center} 
\mbox{\epsfig{file=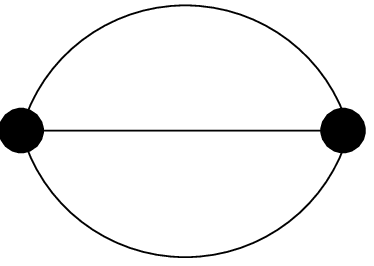,width=1.5cm,height=3cm,angle=270}}
\label{graph10}
\end{center}
\end{figure}

The stabilizer is the whole group $G_7$; it is clear that there exists only
such a curve, so that the real contribution is $1$; our formula gives $\frac
1{12}$.

\begin{figure}[h]
\begin{center} 
\mbox{\epsfig{file=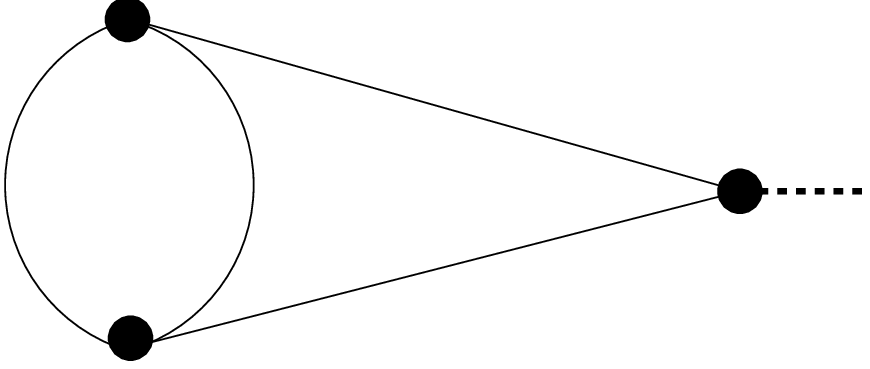,width=2cm,height=4cm,angle=270}}
\label{graph11}
\end{center}
\end{figure}

The stabilizer is $H_7\times K_7^{\prime }$, where $K_7^{\prime }$ is a
subgroup of order $2$ of $K_7$; the contribution is 
\begin{eqnarray*}
\sum_{n\geq 1}\chi \left( \frac{{\cal M}_{0,n+2}}{S_2}\right) \frac{D^n}{n!}
&=&D+\frac E2-\frac D2+\frac{D^2}4 \\
&=&\frac E2+\frac D2+\frac{D^2}4
\end{eqnarray*}
and the formula gives $\frac E4$.

\begin{figure}[h]
\begin{center} 
\mbox{\epsfig{file=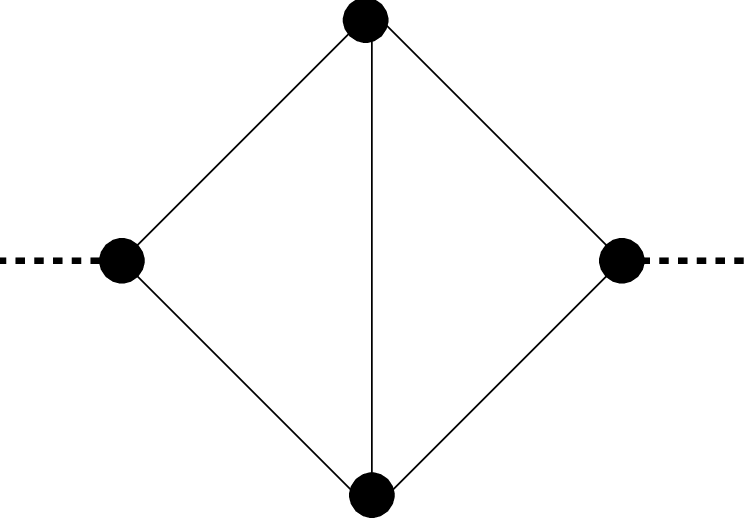,width=1.3cm,height=3cm,angle=270}}
\label{graph12}
\end{center}
\end{figure}

The stabilizer is $H_7$, and the contribution is 
\begin{equation*}
\frac 12\sum_{1\leq n,m\leq 2}\chi \left( \frac{{\cal M}_{0,n+2}\times {\cal %
M}_{0,m+2}}{S_2}\right) \frac{D^{n+m}}{n!m!}=\frac{D^2}2+\frac{D^4}8\text{,}
\end{equation*}
where in the formula we get $\frac 14\left( D^2-D^3+\frac{D^4}4\right) =%
\frac{D^2}4-\frac{D^3}4+\frac{D^4}{16}$.

\begin{figure}[h]
\begin{center} 
\mbox{\epsfig{file=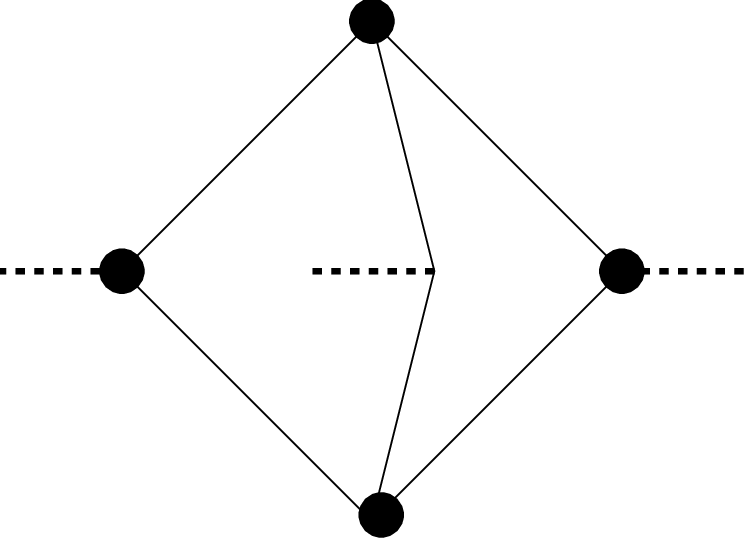,width=1.3cm,height=3cm,angle=270}}
\label{graph13}
\end{center}
\end{figure}

\ Still $H_7$ stabilizes all the graphs of this kind, and the real
contribution is the following: 
\begin{equation*}
\frac 16\sum_{1\leq n,m,p\leq 2}\chi \left( \frac{{\cal M}_{0,n+2}\times 
{\cal M}_{0,m+2}\times {\cal M}_{0,p+2}}{S_2}\right) \frac{D^{n+m+p}}{n!m!p!}%
=\frac{D^3}6+\frac{D^5}8
\end{equation*}
The formula gives instead $\frac 1{12}\left( D^3-\frac 32D^4+\frac 34D^5-%
\frac{D^6}8\right) =\frac{D^3}{12}-\frac{D^4}8+\frac{D^5}{16}-\frac{D^6}{96}$.

\begin{figure}[h]
\begin{center} 
\mbox{\epsfig{file=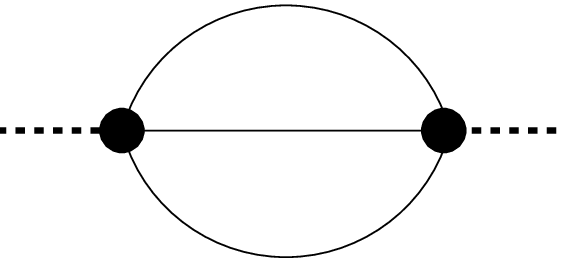,width=1.3cm,height=3cm,angle=270}}
\label{graph14}
\end{center}
\end{figure}

(we assume that at least one vertex has valence $\geq 4$). The stabilizer is 
$K_7$, and the contribution is 
\begin{equation*}
\frac 12\sum\begin{Sb} 0\leq n,m\leq 2  \\ \left( n,m\right) \neq \left(
0,0\right)  \end{Sb}  \chi \left( \frac{{\cal M}_{0,n+3}\times {\cal M}%
_{0,m+3}}{S_3}\right) \frac{D^{n+m}}{n!m!}=D+\frac 32D^2+\frac{D^3}2+\frac{%
D^4}4\text{,}
\end{equation*}
whereas the formula gives $\frac 1{12}\left( -2D+3D^2-2D^3+D^4\right)
=-\frac D6+\frac{D^2}4-\frac{D^3}6+\frac{D^4}{12}$.

\begin{figure}[h]
\begin{center} 
\mbox{\epsfig{file=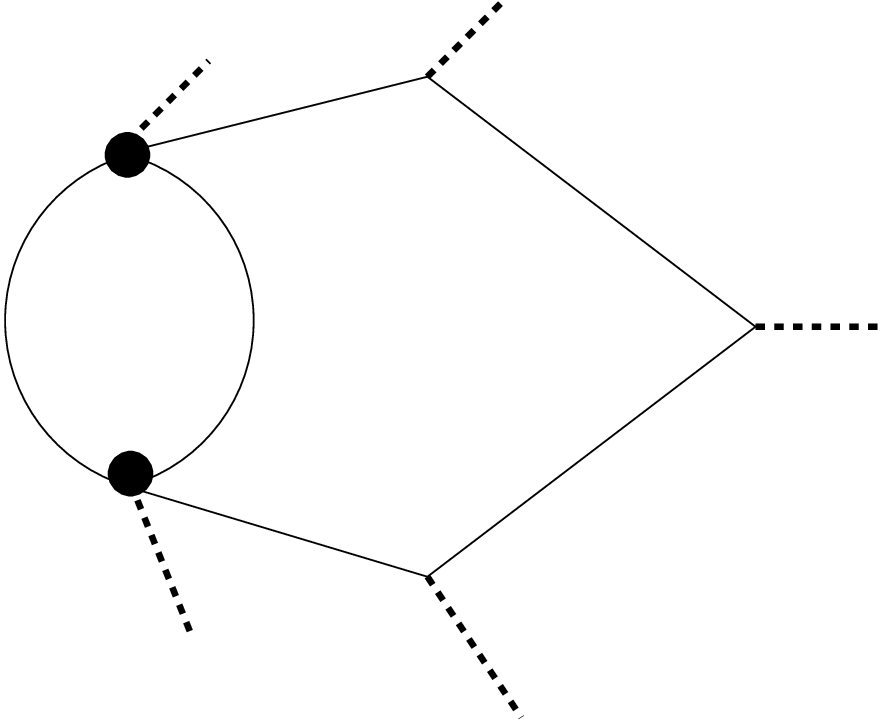,width=3cm,height=5cm,angle=270}}
\label{graph15}
\end{center}
\end{figure}

The stabilizer is a subgroup of order $2$ of $K_7$, and we should put 
\begin{equation*}
\frac 12\left( \sum\begin{Sb} 0\leq n,m\leq 1  \\ \left( n,m\right) \neq
\left( 0,0\right)  \end{Sb}  \chi \left( \frac{{\cal M}_{0,n+3}\times {\cal M%
}_{0,m+3}}{S_2}\right) \frac{D^{n+m}}{n!m!}\right) ED^{\prime }=\frac{ED^2}{%
2\left( 1-E\right) }
\end{equation*}
instead of $\frac{E\left( -2D+D^2\right) }{4\left( 1-E\right) }$, and, for
the case $n=m=0$, since we already corrected part of this term, we
substitute $\frac{E^2}{2\left( 1-E\right) }$ to $\frac{E^2}{4\left(
1-E\right) }$.

The final contribution of type $7$ graphs is 
\begin{eqnarray*}
&&\ \frac 1{12}\frac 1{\left( 1+D\right) ^2}\frac 1{\left( 1-E\right) ^3}+%
\frac{11}{12}+\frac E2+\frac D2+\frac{D^2}4-\frac E4 \\
&&+\frac{D^2}2+\frac{D^4}8-\frac{D^2}4+\frac{D^3}4-\frac{D^4}{16} \\
&&+\frac{D^3}6+\frac{D^5}8-\frac{D^3}{12}+\frac{D^4}8-\frac{D^5}{16}+\frac{%
D^6}{96} \\
&&+D+\frac 32D^2+\frac{D^3}2+\frac{D^4}4+\frac D6-\frac{D^2}4+\frac{D^3}6-%
\frac{D^4}{12} \\
&&+\frac{ED^2}{2\left( 1-E\right) }\ -\frac{-2DE+ED^2}{4\left( 1-E\right) }+%
\frac{E^2}{4\left( 1-E\right) } \\
&=&\frac 1{12\left( 1+D\right) ^2\left( 1-E\right) ^3}+\frac{E^2}{4\left(
1-E\right) }+\frac{E\left( 2D+D^2\right) }{4\left( 1-E\right) } \\
&&+\frac{11}{12}+\frac 53D+\frac 78D^2+D^3+\frac{17}{48}D^4+\frac{D^5}{16}+%
\frac{D^6}{96}\text{.}
\end{eqnarray*}
\bigskip

Now we are able to write down the complete generating function for genus $2$%
: 
\begin{eqnarray*}
 K_2 &=&\frac{1}{1440(1+D)^2(E-1)^3}[ -2D^8(E-1)^2(7+3E)-24D^7(E-1)^2(-7+17E) \\
&&+30D^5\left( E-1\right) ^2\left( 61E-221\right) -3D^6\left( E-1\right)
^2\left( 259+201E\right)  \\
&&+360D\left( 45E^3-167E^2+206E-84\right) +60\left(
73E^3-270E^2+336E-144\right)  \\
&&+180D^2\left( 138E^3-519E^2+635E-254\right) +60D^3\left(
341E^3-1322E^2+1633E-652\right)  \\
&&+15D^4\left( 631E^3-2640E^2+3395E-1386\right) ]\text{.}
\end{eqnarray*}

By developing in power series, we get: 
\begin{equation*}
K_2\left( t\right) =6+13t+21t^2+\frac{181}6t^3+\frac{251}{6}t^4+\frac{6853%
}{120}t^5+\frac{27971}{360}t^6+\frac{177673}{1680}t^7+o(t^8)\text{;}
\end{equation*}
from this we read 
\begin{equation*}
\begin{tabular}{|lllllllll|}
\hline
\multicolumn{1}{|l|}{$n$} & \multicolumn{1}{l|}{$0$} & \multicolumn{1}{l|}{$%
1 $} & \multicolumn{1}{l|}{$2$} & \multicolumn{1}{l|}{$3$} & 
\multicolumn{1}{l|}{$4$} & \multicolumn{1}{l|}{$5$} & \multicolumn{1}{l|}{$6$%
} & $7$ \\ \hline
\multicolumn{1}{|l|}{$\chi \left( \overline{{\cal M}}_{2,n}\right) $} & 
\multicolumn{1}{l|}{$6$} & \multicolumn{1}{l|}{$13$} & \multicolumn{1}{l|}{$%
42$} & \multicolumn{1}{l|}{$181$} & \multicolumn{1}{l|}{$1004$} & 
\multicolumn{1}{l|}{$6853$} & \multicolumn{1}{l|}{$55942$} & $533019$ \\ 
\hline
\end{tabular}
\end{equation*}


\begin{thebibliography}{AC}
\bibitem[AC]{AC}  E. Arbarello, M. Cornalba, {\it Calculating cohomology
groups of moduli spaces of curves via algebraic geometry}, math.AG 9803001,
preprint SNS (1998),

\bibitem[FH]{FH} W. Fulton, J. Harris, {\it Representation Theory }, GTM 129, Springer (1991),

\bibitem[Ga]{Ga} G. Gaiffi, {\it The actions of  $S_{n+1}$ and $S_n$  on the
cohomology ring of a Coxeter arrangement of type $A_{n-1}$ },
Manuscripta Mathematica  91 (1996), 83-94,

\bibitem[G1]{G1}  E. Getzler, {\it The semi-classical approximation for
modular operads}, math.AG 9612005, preprint (1996),

\bibitem[G2]{G2}  E. Getzler, {\it Topological recursion relations in genus 2}, math.AG 9801003,
preprint (1998),

\bibitem[HZ]{HZ}  J. Harer, D. Zagier, {\it The Euler characteristic of the
moduli space of curves}, Invent. math. 85, 457-485 (1986),

\bibitem[L]{L} G. I. Lehrer, {\it On the Poincar\'{e} series associated with Coxeter group
actions on complements of hyperplanes}, J. London Math. Soc.(2) 36 (1987),
275-294,

\bibitem[Ma]{Ma} O. Mathieu, {\it Hidden $\Sigma_{n+1}$-actions}, Commun. Math. Phys. 176 (1996), 467-474.


\end{thebibliography}
\end{document}